\documentclass[a4paper,11pt]{article}

\usepackage{amsmath}
\usepackage{amsfonts}
\usepackage{amssymb}
\usepackage{amsthm}
\usepackage{graphicx}
\usepackage{verbatim}
\usepackage{color}
\usepackage{authblk}
\usepackage{setspace}

\newtheorem{thm}{Theorem}[section]

\newtheorem{conj}{Conjecture}[part]

\theoremstyle{definition}
\newtheorem{defin}[thm]{Definition}

\newcommand{\eps}{\varepsilon}
\newcommand{\B}{\mathbb{B}}
\newcommand{\C}{\mathbb{C}}
\newcommand{\D}{\mathbb{D}}
\newcommand{\N}{\mathbb{N}}

\newcommand{\R}{\mathbb{R}}
\newcommand{\T}{\mathbb{T}}
\newcommand{\W}{\mathbb{W}}
\newcommand{\Z}{\mathbb{Z}}
\newcommand{\BB}{\mathcal{B}}
\newcommand{\DD}{\mathcal{D}}

\newcommand{\RR}{\mathcal{R}}
\newcommand{\RHH}{\mathcal{RH}}
\newcommand{\MM}{\mathcal{M}}
\newcommand{\XX}{\mathcal{X}}
\newcommand{\TT}{\mathcal{T}}

\newcommand{\LL}{\mathcal{L}}

\newcommand{\deltabf}{\text{\boldmath $\delta$}}

\author[1]{Pau Rabassa}
\author[2]{Angel Jorba}
\author[2]{Joan Carles Tatjer}
\affil[1] {\mbox{Johann Bernoulli Institute for Mathematics and Computer Science,}
\centerline{ \mbox{University of Groningen, Groningen, The Netherlands}}
\newline  
\mbox{E-mail: {\tt paurabassa@gmail.com}}                                      
\vspace{2mm}}
\affil[2]{\mbox{Departament of Matem\`atica Aplicada i An\`alisi,}
\mbox{Universitat de Barcelona, Barcelona, Spain} 
\mbox{E-mails: {\tt angel@maia.ub.edu}, {\tt jcarles@maia.ub.es}}}

\date{}

\title{
Towards a renormalization theory for quasi-periodically forced 
                one dimensional maps {III}. {N}umerical 
 		Support\thanks{This work has been supported by the 
		MEC grant MTM2009-09723 and the CIRIT grant 2009 SGR 67.
                P.R. has been partially supported by the PREDEX project, funded
                by the Complexity-NET: {\tt www.complexitynet.eu}.}
}
\begin{document}

\maketitle
\begin{abstract}
In a previous work by the authors the one dimensional 
(doubling) renormalization operator was extended to the case of quasi-periodically 
forced one dimensional maps. The theory was used to explain different 
self-similarity and universality observed numerically in the 
parameter space of the Forced Logistic Maps. The extension proposed 
was not complete in the sense that we assumed a total of four conjectures 
to be true.  In this paper we present numerical support for these conjectures. 
We also discuss the applicability of this theory to the Forced Logistic Map. 
\end{abstract}

\tableofcontents

\section{Introduction}

This is the third of a series of papers (together with
\cite{JRT11a} and \cite{JRT11b}) proposing an extension 
of the one dimensional renormalization theory for the 
case of quasi-periodically forced one dimensional maps.
These three papers are closely related but each of them has been 
written to be readable independently. See also \cite{Rab10} for 
a more detailed discussion. 
In the previous two papers we were concerned with the theoretical 
part of the theory. In this paper we include different 
numerical computations which support the conjectures 
introduced for the developing of this theory. To do that 
we briefly review the theory developed in the previous two papers, 
skipping technical details and we adding some numerical computations 
to the discussion.

The universality and self-renormalization properties in the 
cascade of period doubling bifurcations of families of unimodal 
maps is a well known phenomenon. The paradigmatic example  
is the Logistic map $l_\alpha(x) = \alpha x(1-x)$. 
Given a typical one parametric family 
of unimodal maps $\{f_\alpha\}_{\alpha\in I}$  one  observes
numerically that there exists a sequence of parameter
values $\{d_n\}_{n\in \N} \subset I$ such that
the attracting periodic orbit of the map undergoes a period doubling
bifurcation. Between one period doubling and the next one there
exists also a parameter value $s_n$, for which the critical point
of $f_{s_n}$ is a periodic orbit with period $2^n$. One can
observe numerically that
\begin{equation}
\label{universal limit sumicon}
\lim_{n\rightarrow \infty} \frac{d_n - d_{n-1}} {d_{n+1} - d_{n}} = 
\lim_{n\rightarrow \infty} \frac{s_n - s_{n-1}} {s_{n+1} - s_{n}} =
\deltabf  = \texttt{ 4.66920...}. 
\end{equation}
This convergence indicates a self-similarity 
on the parameter space of the family. On the other hand, 
the constant $\deltabf$ is universal,
in the sense that for any family  of unimodal maps with a quadratic
turning point having a cascade of period doubling bifurcations,
one obtains the same ratio $\deltabf$.

To explain these phenomena Collet and Treser (\cite{CT78})
and  Feigenbaum (\cite{Fei78,Fei79}) proposed simultaneously
the renormalization operator. Their explanation was based on the 
existence of a hyperbolic fixed point of the operator with suitable 
properties. The first proof of the existence of this 
point and its hyperbolicity were obtained with numerical
assistance \cite{Lan82,EW87}. A decade later, Sullivan (see \cite{Sul92})
generalized the operator and provided a theoretical 
proof of the hyperbolicity using complex dynamics techniques. See 
\cite{Lyu99, MvS93} for extensive summaries on the theory.

In \cite{JRT11p2} we presented numerical evidences of self similarity 
and universality for families of quasi-periodically forced Logistic maps. 
These are maps in the cylinder where the dynamics 
on the periodic variable are given by a rigid rotation and 
the dynamics on the other variable are given by  the Logistic Map
plus a small perturbation which depends 
on both variables. This kind of maps have its origins 
in studies related to the existence of strange non-chaotic attractors
(see \cite{Bje09, FKP06, FH10, HH94, Jag03, PMR98}). Let us 
describe these numerical evidences with more detail. 

\subsection{Numerical evidence of self-similarity and universality for
quasi-periodic forced maps}
\label{subsection intro numerical evidences}

\begin{figure}[t]
\centering
\resizebox{16cm}{!}{
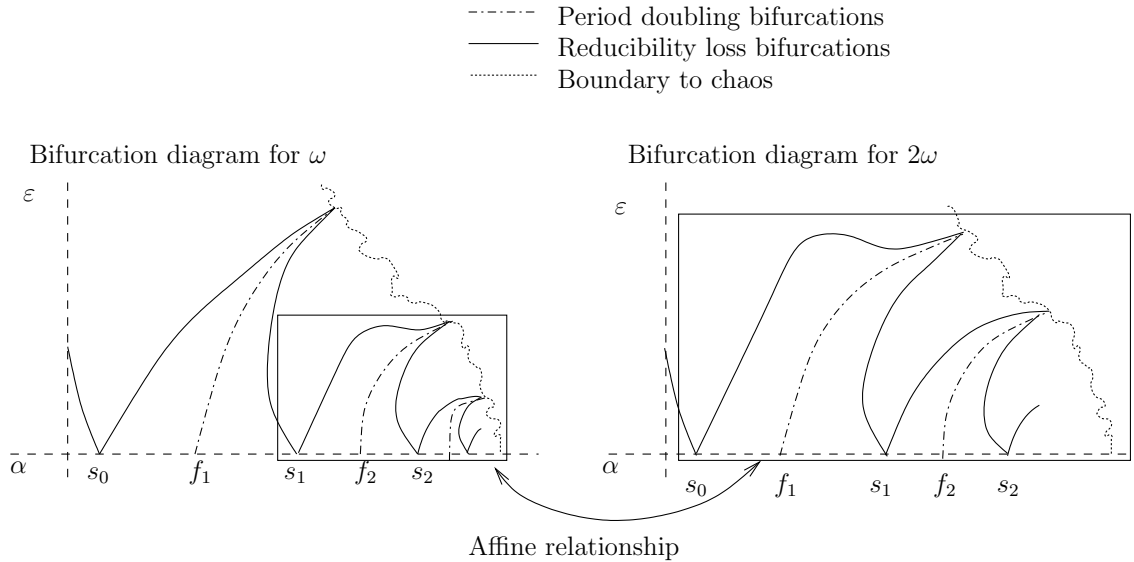
\phantom{aaaaa}
}
\caption{
Schematic representation of the bifurcations diagram
of the Forced Logistic Map, for rotation number equal to
$\omega$ (left) and $2\omega$ (right). See the text for more
details.}
\label{Esquema BifPar}
\end{figure}

Consider $\{g_{\alpha,\eps}\}_{(\alpha,\eps) \in J\subset\R^2}$ 
a two parametric family of quasi-periodic maps in the cylinder 
$\T\times \R$ of the form
\begin{equation}
\label{q.p- family}
\left.
\begin{array}{rcl}
\bar{\theta} & = & \theta + \omega  ,\\
\bar{x} & = & \alpha x(1-x) +  \eps h(\theta,x)  ,
\end{array}
\right\}
\end{equation}
with $\omega$ a Diophantine number, $\alpha$ and $\eps$ parameters and 
$h$ a periodic function with respect to $\theta$ which can also depend 
on $\alpha$ and $\eps$.  Recall that the Logistic map $\bar{x}=\alpha x(1-x)$ has  
a complete cascade of period doubling bifurcations. As before, let
 $\{d_n\}_{n\in \N} \subset I$ denote the parameter values where
the attracting periodic orbit undergoes a period doubling
bifurcation and $\{s_n\}_{n\in \N} \subset I$ the values for which the critical point
of $f_{s_n}$ is a periodic orbit with period $2^n$.

In \cite{JRT11p} we computed some bifurcation diagrams in terms
of the dynamics of the attracting set. We have  taken
into account different properties of the attracting set, as the
Lyapunov exponent  and, in the case of having
a periodic invariant curve, its period and its reducibility. 
The reducibility loss of an invariant curve is not a bifurcation 
in the classical sense that the attracting set of the map changes 
dramatically, only the spectral properties of the transfer operator 
associated to the continuation of that curve does (see \cite{JT08}). 
Despite of this, it can be characterized as a bifurcation
(see definition 2.3 in \cite{JRT11p}). 
The numerical computations in \cite{JRT11p} reveal that the parameter
values for which the invariant curve doubles its period are
contained in regions of the parameter space where the invariant
curve is reducible. These computations also reveal that from every 
parameter value $(\alpha,\eps) = (s_n, 0)$ two curves of reducibility 
loss (of the $2^n$-periodic invariant curve) are  born.
This situation is sketched in the 
left panel of figure \ref{Esquema BifPar}.

Assume that these two curves can be locally expressed as
$(s_n + \alpha_n'(\omega)\eps + O(\eps^2),\eps)$  and
$(s_n + \beta_n'(\omega) \eps + O(\eps^2),\eps)$. In \cite{JRT11a} we proved that 
these curves really exist under suitable hypothesis. We have 
also given explicit expressions of the slopes $\alpha_n'(\omega)$ and 
$\beta_n'(\omega)$ in terms of the quasi-periodic renormalization 
operator (introduced there). 
We only focus on $\alpha_n'(\omega)$, but the discussion for $\beta_n'(\omega)$
is completely analogous.

The slopes $\alpha_n'(\omega)$ can be used for the numerical detection 
of universality and self-renormalization phenomena. If the bifurcation 
diagram is self-similar by an affine ratio one should have that $\alpha_n'(\omega)/ 
\alpha_{n-1}'(\omega)$ converges to a constant. In \cite{JRT11p2} 
we compute numerically this ratios and we show that this is not true  
due to the fact that when the period is doubled, the rotation number 
of the system also is. What we find is that there exists an affine relationship 
between the bifurcation diagram of the family for rotation number $\omega$
and the bifurcation diagram of the same family for rotation number
 $2\omega$. This is sketched in figure \ref{Esquema BifPar}.

Concretely, in \cite{JRT11p2} we observed numerically the following behavior. 

\begin{itemize}
\item {\bf First numerical observation:}  the sequence $\alpha_n'(\omega)/ 
\alpha_{n-1}'(\omega)$ is not convergent in $n$. But, for $\omega$ fix,
one obtains the same sequence for any family of quasi-periodic forced map
like (\ref{q.p- family}), with
a quasi-periodic forcing of the type $h(\theta,x) =  f_1(x) \cos(\theta)$.
\item {\bf Second numerical observation:}  the sequence $\alpha_n'(\omega)/ 
\alpha_{n-1}'(2\omega)$ associated to maps like (\ref{q.p- family})
is convergent in $n$ when we take the quasi-periodic forcing of the 
type $h(\theta,x) = f_1(x) \cos(\theta)$. 
The limit depends on $\omega$ and $f_1$. 
\item {\bf Third numerical observation:} the two previous observations
are not true when the quasi-periodic forcing is of the
type $h_\eta(\theta,x) = f_1(x) \cos(\theta) + \eta f_2(x) \cos(2\theta)$
when $\eta\neq 0$. But the sequence $\alpha_n'(\omega)/ 
\alpha_{n-1}'(2\omega)$ associated to the map (\ref{q.p- family})
with $h = h_\eta$ is $\eta$-close to the same maps with $h= h_0 $
\end{itemize}

In \cite{JRT11a} we extended the renormalization operator and we obtained 
explicit expressions of the slopes $\alpha_n'(\omega)$ and $\beta_n'(\omega)$ 
in terms of the quasi-periodic renormalization operator. This is 
reviewed in section \ref{section review q-p renor}. In \cite{JRT11b} 
we give a theoretical 
explanation to the numerical observations described above in terms of
the dynamics of the quasi-periodic  renormalization operator.
This is reviewed in section \ref{section A not yet rigorous explanation}. 
The novelty in this paper is that we present numerical support to the 
conjectures done in \cite{JRT11a, JRT11b}. This numerical support is presented 
after the statement of each of the conjectures, with the exception 
of conjecture {\bf \ref{conjecture H2}}, which is given in section 
\ref{section applicability FLM}. 
In Appendix \ref{section numerical computation of the specturm of L_omega}
we describe the numerical approximation used to discretize the 
renormalization operator and how we use it to compute the spectrum of its 
derivative.

\section{{E}xistence of reducibility loss bifurcations}
\label{section review q-p renor} 


Consider a quasi-periodic forced map like 
\begin{equation}
\label{q.p. forced system interval}
\begin{array}{rccc}
F:& \T\times I &\rightarrow & \T \times I \\
  & \left( \begin{array}{c} \theta \\ x \end{array}\right)
  & \mapsto
  & \left( \begin{array}{c} \theta + \omega \\ f(\theta,x) \end{array}\right),
\end{array} 
\end{equation}
with $f\in C^r(\T\times I ,I)$. To define the renormalization 
operator it is only necessary that $r \geq 1$, but 
we restrict this study to the analytic case due to technical reasons.
Along section \ref{section definition and basic properties} 
it is not necessary to require $\omega$  Diophantine, but 
it will be necessary in section \ref{chapter application}. 

The definition of the operator is done in a perturbative way, 
in the sense that it is only applicable to maps 
$f(\theta,x) = g(x) + h(\theta,x)$ with $g$ renormalizable in 
the one dimensional sense and $h$ small. 

\subsection{Definition of the operator and basic properties} 
\label{section definition and basic properties}

\subsubsection{Preliminary notation} 

Let $\W$ be an open set in the complex plane containing the
interval $I_\delta=[-1-\delta, 1+\delta]$ and let
$\B_\rho = \{z = x + i y\in \C \text{ such that } |y| < \rho\}$.   Then 
consider  $\BB=\BB(\B_\rho,\W)$ the space of functions
$f: \B_\rho \times  \W \rightarrow \C$ such that:
\begin{enumerate}
\item  $f$ is holomorphic in $\B_{\rho}\times \W$ and continuous
in the closure of  $\B_{\rho}\times \W$.
\item  $f$ is real analytic.
\item  $f$ is $1$-periodic in the first variable.
\end{enumerate}
This space endowed with the supremum norm is Banach. 

Let $\RHH(\W)$  denote the space of functions
$f:\W \rightarrow \C$ such that are holomorphic
in $\W$, continuous in the closure of  $\W$,
and send real number to real numbers. This space is also 
Banach with the supremum norm. 

Consider the operator
\begin{equation}
\label{equation projection 0}
\begin{array}{rccc}
p_0:& \BB &\rightarrow & \RHH(\W) \\
\displaystyle \rule{0pt}{3ex} & f(\theta,x) & \mapsto &  \rule{0pt}{3ex}
\displaystyle  \int_{0}^{1} f(\theta, x) d\theta .
\end{array} 
\end{equation}
Let $\BB_0$ the natural inclusion of $\RHH(\W) $ 
into $\BB$. Then we have that $p_0$  as a map from $\BB$ to $\BB_0$ 
is a projection  ($(p_0)^2 = p_0$).

\subsubsection{Set up of the one dimensional renormalization operator. }

First we give a concrete definition of the 
one dimensional renormalization operator before extending 
it to the quasi-periodic case. Actually, this  is a minor 
modification of the one given in \cite{Lan82}. 

Given a small value $\delta>0$, let  $\MM_\delta$ denote the 
subspace of $\RHH(\W)$ formed by the even functions $\psi$ which 
send the interval $I_\delta=[-1-\delta,1+\delta]$ into itself, and such that
$\psi(0)=1$ and  $x \psi'(x) <0$ for $x\neq 0$.

Set $a=\psi(1)$, $a'= (1+\delta)a$ and $b'=\psi(a')$. We 
can define $\DD(\RR_\delta)$ as the set of $\psi \in \MM_\delta$ 
such that $a<0$, 
$1> b'>-a'$,
and  $\psi(b') <- a'$.

We define the renormalization operator, $\RR_\delta : 
\DD(\RR_\delta) \rightarrow \MM_\delta$ as 
\begin{equation}
\label{renormalization operator lanford}
\RR_\delta(\psi) (x) =  \frac{1}{a} \psi \circ \psi (a x).
\end{equation}
where $a=\psi(1)$.

Note that, given $\psi \in \DD(\RR_\delta)$, one needs to ensure that 
$\psi\left(a\W\right)\subset \W$ in order to have $\RR_\delta(\psi)$ well defined.
With this aim, let us consider the following hypothesis. 

\begin{description}
\item[H0)] There exists an open set $\W\subset \C$ containing 
$I_\delta$ and a function $\Phi\in \BB $
such that $\phi= p_0(\Phi) $ is a fixed point of the  
renormalization operator $\RR_\delta$  and such that 
the closure of both $a\W$ and $\phi(\Phi)(a \W)$ is 
contained in $\W$ (with $a:=\Phi(1)$). 
\end{description}

In \cite{Lan92}, Lanford claims that the hypothesis {\bf H0} is satisfied 
by the set
\begin{equation}
\label{set Lanford}
\left\{ z\in \C \text { such that } |z^2 -1| < \frac{5}{2} \right\}.
\end{equation}
This set is convenient for him because he works in the
set of even holomorphic functions. 

In \cite{Lan82} Lanford introduces a discretization of the 
(one dimensional) renormalization operator to give
a computer assisted proof of the contractivity of the operator. 
In the present paper we use the same techniques to discretize the 
quasi-periodic renormalization operator, although we do it 
without the use of rigorous interval arithmetics. More details on this
are given in the Appendix \ref{section numerical computation of the specturm of L_omega}. 
We can use this discretization to check the hypothesis {\bf H0} 
for a suitable set $\W$. 

\begin{figure}[t]
\begin{center}
\includegraphics[width=7.5cm]{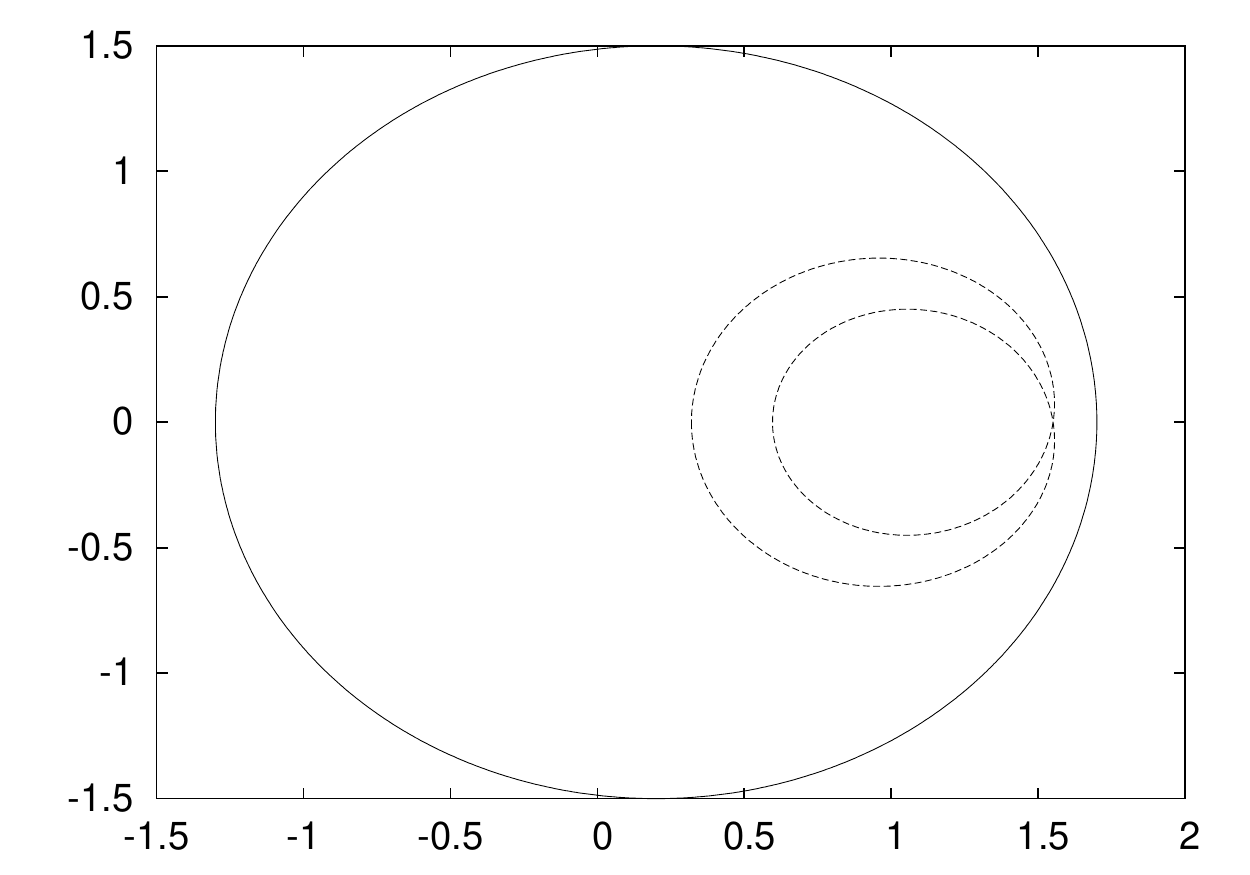}
\caption{ 
Boundaries of the sets $\D\left(\frac{1}{5},\frac{3}{2} \right)$
(solid line)  and   $\phi\left(a \D\left(\frac{1}{5},\frac{3}{2} \right)\right)$
(dashed line).}
\label{figure inclusion boundaries}
\end{center}
\end{figure}

Our study is not restricted to the case of even functions, 
therefore the set (\ref{set Lanford}) used 
by Lanford is not valid in our case. 
Using the method described in the Appendix 
\ref{section numerical computation of the specturm of L_omega} 
we recomputed the fixed point $\phi$
of the (one dimensional) renormalization operator $\RR$. 
The fixed point has been computed by means of a Newton
method with our discretization and then we have
checked that the Taylor expansion around zero coincides 
with the one given in \cite{Lan82}.

With this approximation of $\phi$ we checked (numerically) that 
$\D\left(\frac{1}{5},\frac{3}{2} \right)$ the disc of the complex 
plane  centered at $\frac{1}{5}$ with radius $\frac{3}{2}$ satisfies 
the conditions required to the set $\W$ in hypothesis {\bf H0}. In other
words, we checked that 
$\phi\left( 
a \D\left(\frac{1}{5},\frac{3}{2} \right)\right)$ is
contained inside  $\D\left(\frac{1}{5},\frac{3}{2} \right)$
(recall that $a=\phi(1)$ and $a \D(z_0,\rho) = \{ z \in \C \thinspace 
| \thinspace az \in \D(z_0,\rho) \}$ ).

Denote by $\partial \D(z_0,\rho) $ the boundary of the disk
$\D(z_0,\rho) $. In figure
\ref{figure inclusion boundaries} we have
plotted the sets $\partial \D\left(\frac{1}{5},\frac{3}{2} \right)$
 and   $\phi\left(a \partial \D\left(\frac{1}{5},\frac{3}{2} \right)\right)$
which give a visual evidence of the inclusion. Recall that $\phi$ is
analytic, then to check that the set
$\D\left(\frac{1}{5},\frac{3}{2} \right)$ is mapped inside the
set delimited by $\phi\left(a \partial \D\left(\frac{1}{5},\frac{3}{2} \right)\right)$ it is
enough to check that one point in the interior of
$\D\left(\frac{1}{5},\frac{3}{2} \right)$ is mapped in the
interior of $\phi\left(a \partial \D\left(\frac{1}{5},\frac{3}{2} \right)\right)$.
Recall that $\phi(0)=1$ by hypothesis, therefore the inclusion holds.

\subsubsection{Definition of the renormalization operator for  
quasi-periodically forced maps}

Consider the space $\XX\subset \BB$ defined as: 
\[
\XX= \{ f \in C^r(\T \times I_\delta, I_\delta) | \thinspace 
 p_0(f) \in \MM_\delta\}. 
\]
Consider also the decomposition $\XX=\XX_0 \oplus \XX_0^c$ given by the projection 
$p_0$. In other words, we have $\XX_0=\{f \in \XX \thinspace |\text{ } p_0(f)=f\}$ and 
$\XX^c_0=\{f \in \XX \thinspace | \text{ } p_0(f)=0\}$. Note 
that from the definition of $\XX$ it follows that $\XX_0$ is an isomorphic 
copy of $\MM_\delta$.  

Given a function $g\in \XX$, we  
define 
the {\bf quasi-periodic renormalization} of $g$ as 
\begin{equation}
\label{operator tau}
[\TT_\omega(g)](\theta,x) :=  \frac{1}{\hat{a}} g(\theta + \omega,
g(\theta, \hat{a}x)),
\end{equation}
where $\displaystyle \hat{a} = \int_{0}^{1} 
g(\theta, 1) d\theta$. 

Then we have that there exists a set $\DD(\TT)$, 
open in $\left(p_0 \circ \TT_\omega\right)^{-1} (\MM_\delta)$, 
where the operator is well defined. 
Moreover this set contains  $\DD_0(\TT)$, the inclusion of $\DD(\RR)$ in 
$\BB$. By definition we have that $\TT_\omega$ restricted 
to $\DD_0(\TT)$ is isomorphically conjugate to $\RR$, therefore 
the fixed points of $\RR$ extend to fixed points of $\TT_\omega$. 
Assume that {\bf H0} holds and let $\Phi$ be 
the fixed point given by this hypothesis. Then we have that there exists 
$U\subset \DD(\TT)\cap\BB$, an open neighborhood of $\Phi$, 
such that $\TT_\omega :U\rightarrow \BB$ is well defined. 
Moreover we have that $\TT_\omega$ is Fr\'echet differentiable for any  $\Psi \in U$.

\subsubsection{Fourier expansion of $D\TT_\omega(\Psi)$.} 
\label{section The Fourier expansion of DT} 

Let $\Psi$ be a function in a neighborhood of $\Phi$ (given in hypothesis
{\bf H0}) where $\TT_\omega$ is differentiable. Additionally, 
assume that $\Psi\in \DD_0(\TT_\omega)$.

Given a function $f\in \BB$ we can consider its complex Fourier
expansion in the periodic variable
\begin{equation}
\label{fourier expansion complex}
f(\theta,z)= \sum_{k\in \Z} c_k(z) e^{2\pi k\theta i }, 
\end{equation}
with
\[
c_k(z)= \int_{0}^{1} f(\theta,z) e^{-2\pi k\theta i} d\theta.
\]

Then we have that $D\TT_\omega$ ``diagonalizes''
with respect to the complex Fourier expansion, in the
sense that we have 
\begin{equation}
\label{equation Fourier expansion differential renormalization operator}
\left[D\TT_\omega(\Psi) f\right](\theta,z) = D\RR_\delta[c_0](z) 
+  \sum_{k\in \Z\setminus\{0\}} \left([L_1(c_k)](z) + 
[L_2(c_k)](z) e^{2\pi k\omega i}\right) e^{2\pi k\theta i },
\end{equation}
where
\[
\begin{array}{rccc}
L_1: & \RHH(\W) & \rightarrow & \RHH(\W)  \\
  &   g(z)  & \mapsto & 
\displaystyle  \frac{1}{a} \psi'\circ\psi(a z) g(az),
\end{array} 
\]
and 
\[
\begin{array}{rccc}
L_2: & \RHH(\W) & \rightarrow & \RHH(\W)  \\   &   g(z)  & \mapsto &
\displaystyle  \frac{1}{a} g\circ\psi(a z),
\end{array}
\] with $\psi=p_0(\Psi)$ and 
$a=\psi(1)$. 

An immediate consequence of this diagonalization is the following. Consider 
\begin{equation}
\label{equation spaces bbk}
\BB_k:= \big\{ f \in B |\text{ } f(\theta, x) = u(x) \cos(2\pi k \theta) + 
v(x) \sin(2\pi k\theta), \text{ for 
some } u,v\in \RHH(\W)\big\},
\end{equation} 
then we have that the spaces $\BB_k$ are
invariant by $D\TT(\Psi)$ for any $k>0$. 

Moreover $D\TT_\omega (\Psi)$ restricted to $\BB_k$ is
conjugate to $\LL_{k\omega}$, where $\LL_\alpha$ is the defined 
as 
\begin{equation}
\label{equation maps L_omega}
\begin{array}{rccc}
\LL_\alpha: &  \RHH(\W)\oplus  \RHH(\W) & \rightarrow 
&  \RHH(\W)\oplus  \RHH(\W)  \\ \\
& \left( \begin{array}{c} u \\v \end{array} \right) &
\mapsto &  \left( \begin{array}{c}  L_1(u) \\ L_1(v)
\end{array} \right) +
\left( \begin{array}{cc}  \cos(2\pi \alpha) & - \sin(2\pi \alpha)  \\
\sin(2\pi \alpha) & \cos( 2\pi \alpha) 
\end{array} \right)  \left( \begin{array}{c}  L_2(u) \\
L_2(v) 
\end{array} \right) . 
\end{array}
\end{equation}

Then we have that the understanding of the 
derivative of the renormalization operator in $\BB$ is equivalent 
to the study of the operator $\LL_\omega$ for a any $\omega\in \T$.

\subsubsection{Properties of $\LL_\omega$} 
\begin{figure}[t]
\begin{center}
\includegraphics[width=12cm]{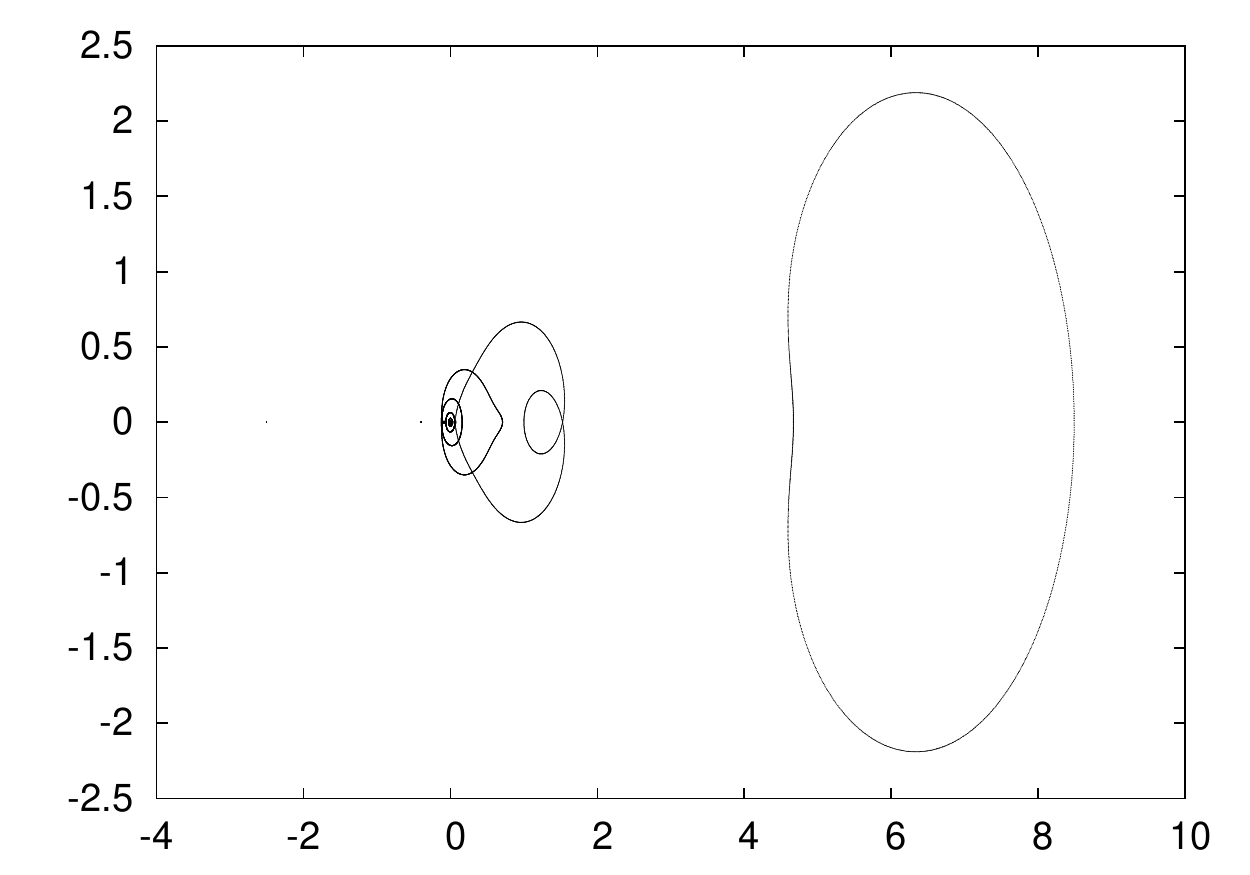}
\includegraphics[width=6cm]{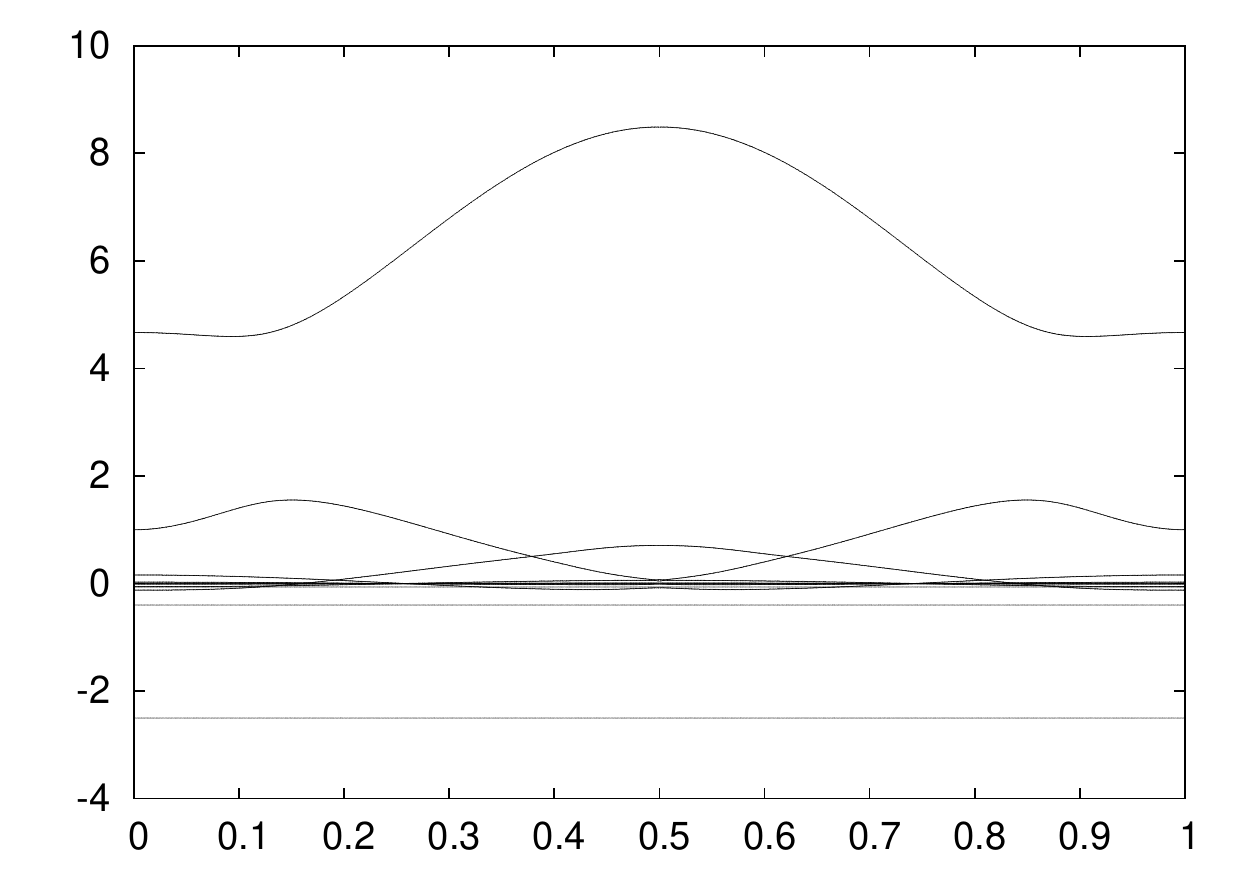}
\includegraphics[width=6cm]{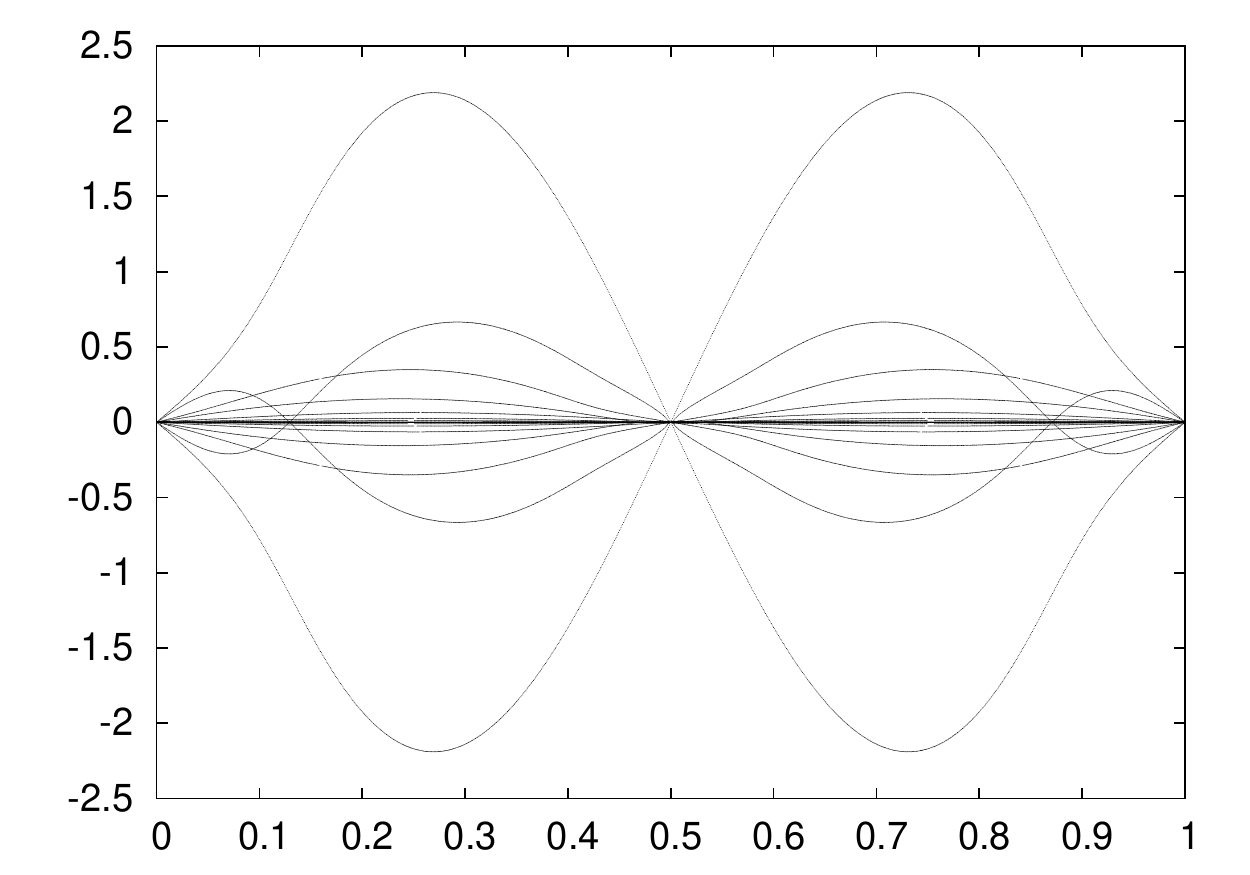}
\caption{
Numerical approximation of the spectrum of $\LL_\omega$
for $\omega \in \T$. Top: projection in the complex plane of the
spectrum when $\omega$ varies in $\T$. Bottom left: evolution of the
real part with respect to $\omega$. Bottom right:  evolution of the
imaginary part with respect to $\omega$.
}
\label{Operator Spectrum}
\end{center}
\end{figure}

Given a value $\gamma \in \T$, consider the rotation $R_{\gamma}$ defined
as
\begin{equation}
\label{equation rotation rgamma}
\begin{array}{rccc}
R_\gamma: &  \RHH(\W)\oplus  \RHH(\W) & \rightarrow
&  \RHH(\W)\oplus  \RHH(\W)  \\ \\
& \left( \begin{array}{c} u \\v \end{array} \right) &
\mapsto &
\left( \begin{array}{cc}  \cos( 2\pi \gamma) & - \sin(2\pi \gamma)  \\
\sin(2\pi \gamma) & \cos(2\pi \gamma)
\end{array} \right)  \left( \begin{array}{c} u \\ v
\end{array} \right) ,
\end{array}
\end{equation}
then we have that $\LL_\omega$ and $R_\gamma$ commute
for any $\omega, \gamma \in \T$. 

This has some consequences on the spectrum of $\LL_\omega$. Concretely, 
we have that any eigenvalue of $\LL_\omega$ (different from zero) is 
either real with geometric multiplicity even, or a pair of complex 
conjugate eigenvalues. On the other hand $\LL_\omega$ depends analytically
on $\omega$, which (using theorems III-6.17 and VII-1.7 of \cite{Kat66})
imply that (as long as the eigenvalues of $\LL_\omega$ are simple) 
the eigenvalues and their associated eigenspaces depend 
analytically on the parameter $\omega$. 

Finally, doing some minor changes on the domain of definition, 
we can prove the compactness of the operator $\LL_\omega$. 
Recall that the compactness of an operator implies that its
spectrum is either finite or countable with $0$ on its 
closure (see for instance theorem III-6.26 of \cite{Kat66}).

In figure \ref{Operator Spectrum} we have a numerical approximation of
the spectrum of the operator $\LL_\omega$ depending on $\omega$. We can 
observe that the properties described above are satisfied. 
The details on the numerical computations involved to approximate
the spectrum are described in Appendix \ref{section numerical computation 
of the specturm of L_omega}. Several numerical tests
on the reliability of the results are also included there.

\subsection{Reducibility loss and quasi-periodic renormalization}
\label{chapter application} 

Given a map $F$ like (\ref{q.p. forced system interval}) with 
$f\in\BB$ and $\omega\in\T$,  we 
denote by $f^{n}:\T \times \R \rightarrow \R$ the $x$-projection
of $F^n(x,\theta)$. Equivalently $f^n$ can be defined through the recurrence
\begin{equation}
\label{definicio falan}
f^n(\theta, x) = f( \theta +(n-1) \omega, f^{n-1}(\theta,x)), \quad 
f^0(\theta,x) \equiv x.
\end{equation}

From this point on, whenever  $\omega$ is used, it is assumed to 
be Diophantine.  Denote by $\Omega= \Omega_{\gamma,\tau}$
the set of Diophantine numbers,
this is the set of $\omega\in \T$
such that there exist $\gamma >0$ and $\tau \geq 1$
such that
\[
|q \omega - p| \geq \frac{\gamma}{|q|^{\tau}}, \quad 
\text{ for all } (p,q) \in \Z \times (\Z \setminus \{0\}). 
\]

Additionally, we will need to assume that the following conjecture is true.  

\begin{conj}
\label{conjecture H2} 
The operator $\TT_{\omega}$ (for any $\omega\in\Omega$) is an injective 
function when restricted to the domain $\BB\cap \DD(\TT)$. Moreover, 
there exists $U$ an open set of $\DD(\TT)$ containing $W^u(\Phi,\RR)
\cup W^s(\Phi,\RR)$\footnote{Here $W^s(\Phi,\RR)$ and $W^u(\Phi,\RR)$ 
are considered as the inclusion in $\BB$ of the 
stable and the unstable manifolds of the fixed point $\Phi$ (given by 
{\bf H0}) by the map $\RR$ in the topology of $\BB_0$ (the inclusion 
of one parametric maps in $\BB$).}
where the operator $\TT_\omega$ is differentiable.
\end{conj}

In \cite{JRT11a} we discuss the difficulties for proving this 
conjecture. A priory there is no way to check numerically this kind 
of conjecture. A posteriori we have that 
the results obtained assuming this conjecture are coherent with the numerical 
computations (see section \ref{section applicability FLM}). 

\subsubsection{Consequences for a two parametric family of maps} 
\label{Section consequences for a two parametric family of maps}

Consider a two parametric families of maps  $\{c(\alpha,\eps)\}_{(\alpha,\eps)\in A}$ 
contained in $\BB$, with $A = [a,b]\times[0,c]$ and $a$, $b$ 
and $c$ are real numbers (with $a<b$ and $0<c$). 
We assume that the dependency on the parameters is analytic.

Consider the following hypothesis on the family of maps. 
\begin{description}
\item[H1)]  The family $\{c(\alpha,\eps)\}_{(\alpha,\eps)\in A}$ uncouples for
$\eps=0$, in the sense that the family $\{c(\alpha,0)\}_{\alpha\in[a,b]}$
does not depend on $\theta$ and it has a full cascade of
period doubling bifurcations. We assume that the family
$\{c(\alpha,0)\}_{\alpha\in[a,b]}$ crosses transversally the
stable manifold of $\Phi$, the fixed point
of the renormalization operator, and each of
the manifolds $\Sigma_n$ for any $n\geq 1$, where
$\Sigma_{n}$ is the inclusion in $\BB$ of the
set of one dimensional unimodal maps with
a super-attracting $2^n$ periodic orbit.
\end{description}

In other words, we assume that the family $c(\alpha,\eps)$ 
can be written as, 
\[
c(\alpha,\eps)=c_0(\alpha) + \eps c_1(\alpha,\eps),
\]
with $\{c_0(\alpha)\}_{\alpha\in[a,b]}\subset \BB_0$  having a full 
cascade of period doubling bifurcations. 

Given a family $\{c(\alpha,\eps)\}_{(\alpha,\eps)\in A}$ satisfying 
the hypothesis {\bf H1},  let $\alpha_n$ be the parameter value for 
which the uncoupled family $\{c(\alpha,0)\}_{\alpha\in[a,b]}$
intersects the manifold $\Sigma_n$. Note that the critical 
point of the map $c(\alpha_n,0)$ is a $2^n$-periodic orbit.
Our main achievement in \cite{JRT11a} is 
to prove that from every 
parameter value $(\alpha_n,0)$ there are born two curves in 
the parameter space, each of them corresponding to a reducibility 
loss bifurcation. If we want to give a more precise statement of 
the result we need now to introduce some technical definitions.

Let $\RHH(\B_\rho,\W)$ denote the space of
periodic real analytic maps from $\B_\rho$ to $\W$ and
continuous in the closure of $\B_\rho$.
Consider a map $f_0\in \BB$ and $\omega \in \Omega$, such that $f$ has a 
periodic invariant curve $x_0$ of rotation number $\omega$
with a  Lyapunov exponent less equal than certain $-K_0<0$. 
Using lemma 3.6 in \cite{JRT11a} we have that 
there exist a neighborhood $V\subset \BB$ of $f_0$ and a map 
$x\in \RHH(\B_\rho, \W)$ such that $x(f)$ is a periodic 
invariant curve of $f$  for any $f\in V$. 
Then we can define
the map $G_1 $ as
\begin{equation}
\label{equation definition G1}
\begin{array}{rccc}
G_1:& \Omega\times V &\rightarrow & \RHH(\B_\rho,\C) \\
\rule{0pt}{3ex} &  (\omega,g) & \mapsto &  
D_x g \big(\theta+\omega,g(\theta, \left[x(\omega,g)\right](\theta))\big) 
D_x g \big(\theta,\left[x(\omega,g)\right](\theta)\big). 
\end{array}
\end{equation}

On the other hand, we can consider the counterpart of the map 
$G_1$ in the uncoupled case. Given a map $f_0\in \BB_0$, consider 
$U\subset \BB_0$ a neighborhood of $f_0$ in the $\BB_0$ topology. 
Assume that $f_0$ has a attracting $2$-periodic 
orbit $x_0\in I$. Let $x=x(f)\in \W$ be the continuation of this 
periodic orbit for any $f\in U$.  We have  that $x$ depends 
analytically on the map, therefore it induces a map $x:U \rightarrow \W$. 
Then if we take $U$ 
small enough we have an analytic map 
$x:U\rightarrow \W$ such that $x[f]$ is a periodic orbit of period 2. 
Now we can consider the map
\begin{equation}
\label{equation definition widehatG1}
\begin{array}{rccc}
\widehat{G}_1:& U\subset \BB_0 &\rightarrow & \C \\
\rule{0pt}{3ex} &  f & \mapsto & D_x f \big(f(x[f])\big) D_x f \big( x[f] \big). 
\end{array}
\end{equation}

Note that $\widehat{G}_1$ corresponds to $G_1$ restricted to 
the space $\BB_0$ (but then $\widehat{G}_1(f)$ has to be 
seen as an element of $\RHH(\B_\rho,\W)$). 

Consider the sequences
\begin{equation}
\label{equation sequences corollary directions}
\begin{array}{rcll} 
\omega_k  & = & 2 \omega_{k-1},  & \text{ for }  k=1,..., n-1. \\ 
\rule{0ex}{4ex} 
f^{(n)}_k & = & \RR\left(f_{k-1}^{(n)}\right),
& \text{ for }  k=1,..., n-1. \\ 
\rule{0ex}{4ex} 
u^{(n)}_k & = & D \RR \left(f^{(n)}_{k-1}\right) u^{(n)}_{k-1}, 
&  \text{ for } k = 1, ..., n-1.\\
\rule{0ex}{4ex} 
v^{(n)}_k & = & D \TT_{\omega_{k-1}}  \left(f^{(n)}_{k-1}\right) v^{(n)}_{k-1}, 
&  \text{ for } k = 1, ..., n-1.
\end{array}
\end{equation}
with
\begin{equation}  
\label{equation sequences directions initial}
f^{(n)}_0  =  c(\alpha_n,0), \quad 
u^{(n)}_0  =  \partial_\alpha c(\alpha_n,0), \quad 
 v^{(n)}_0 =  \partial_\eps c(\alpha_n,0).  
\end{equation}

Note that $f^{(n)}_0= \{c(\alpha,0)\}\cap \Sigma_n$, then 
$f^{(n)}_0$ tends to $W^s(\RR,\Phi)$ when $n$ grow. 
Therefore, the sequence  $\{f^{(n)}_k\}_{0\leq k <n}$ attains to $W^s(\RR,\Phi) 
\cup W^u(\RR,\Phi)$ when $n$ grows and consequently there 
exist $n_0$ such that $\{f^{(n)}_k\}_{0\leq k <n} \subset U$ , 
where $U$ is the neighborhood given in conjecture {\bf \ref{conjecture H2}}. 
If the conjecture is true, then the operator
 $\TT_\omega$ is differentiable in the orbit $\{f^{(n)}_k\}_{0\leq k <n} 
\subset U$.  

Consider the following 
hypothesis. 

\begin{description}
\item[H2)]  The family $\{c(\alpha,\eps)\}_{(\alpha,\eps)\in A}$ is such that 
\[
D G_1 \left(\omega_{n-1}, f^{(n)}_{n-1}\right) 
D\TT_{\omega_{n-2}}\left(f^{(n)}_{n-2}\right) \cdots 
D\TT_{\omega_0}\left(f^{(n)}_0\right) \partial_\eps c(\alpha_n,0),
\] 
has a unique non-degenerate minimum (respectively maximum) as a function from $\T$ to $\R$, 
for any $n\geq n_0$. 
\end{description}

Consider a family of maps $\{c(\alpha,\eps)\}_{(\alpha,\eps)\in A}$
such that the hypotheses {\bf H1} and {\bf H2} are satisfied
and  $\omega_0\in \Omega$. If the conjecture {\bf \ref{conjecture H2}} is
true, then  theorem 3.8 in \cite{JRT11a} asserts that there exists
$n_0$ such that, for any $n\geq n_0$,
there exist two bifurcation curves around the
parameter  value $(\alpha_n, 0)$,  such that they correspond to a
reducibility-loss bifurcation of the $2^n$-periodic invariant curve.
Moreover, these curves are locally expressed as
$(\alpha_n + \alpha_n'(\omega) \eps + o(\eps),\eps)$ and
$(\alpha_n^- + \beta_n'(\omega) \eps + o(\eps), \eps)$ with
\begin{equation}
\label{equation alpha n +}
\alpha'_n (\omega) = 
- \frac{m \left(
DG_1 \left(\omega_{n-1}, f^{(n)}_{n-1}\right) v_{n-1}^{(n)}
\right)
}{\rule{0ex}{3.5ex}  
D \widehat{G}_1 \left(f^{(n)}_{n-1}\right) u_{n-1}^{(n)}} ,
\end{equation}
and
\begin{equation}
\label{equation alpha n -}
\beta'_n(\omega) = 
- \frac{M \left(
D G_1 \left(\omega_{n-1}, f^{(n)}_{n-1}\right) v_{n-1}^{(n)}
\right)
}{\rule{0ex}{3.5ex}  
D \widehat{G}_1 \left(f^{(n)}_{n-1}\right) u_{n-1}^{(n)}} ,
\end{equation}
where $G_1$ and $\widehat{G}_1$ are given by equations
(\ref{equation definition G1}) and (\ref{equation definition widehatG1}),
and $m$ and $M$ are the minimum and the maximum as operators,
that is
\begin{equation}
\label{equation definition of minim} 
\begin{array}{rccc}
m:& \RHH(\B_\rho,\C) &\rightarrow & \R \\
\rule{0pt}{3ex} & g & \mapsto & \displaystyle \min_{\theta \in \T} g(\theta). 
\end{array}
\end{equation}
and
\begin{equation}
\label{equation definition of maximum} 
\begin{array}{rccc}
M:& \RHH(\B_\rho,\C) &\rightarrow & \R \\
\rule{0pt}{3ex} & g & \mapsto & \displaystyle \max_{\theta \in \T} g(\theta). 
\end{array}
\end{equation}

Let us focus again on hypothesis {\bf H2}, which is not intuitive. 
We can introduce a stronger condition which much more easy to 
check. Moreover this conditions is automatically satisfied by maps like 
the Forced Logistic Map. Consider a family of maps 
$\{c(\alpha,\eps)\}_{(\alpha,\eps)\in A}$ as before, 
satisfying hypothesis {\bf H1}. 

\begin{description}
\item[H2')]  The family $\{c(\alpha,\eps)\}_{(\alpha,\eps)\in A}$ is such that
the quasi-periodic perturbation $\partial_\eps c(\alpha,0)$ belongs to the set
$\BB_1$ (see equation (\ref{equation spaces bbk}))  for any value of 
$\alpha$ (with $(\alpha,0)\in A$). 
\end{description}

Then we have that {\bf H2'} implies {\bf H2} (see proposition 
3.10 in \cite{JRT11a}).

\section{{A}symptotic behavior of reducibility loss bifurcations}
\label{section A not yet rigorous explanation} 

Consider a two parametric family of maps $\{c(\alpha,\eps)\}_{(\alpha,\eps)
\in A}$ contained in $\BB$, with $A = [a,b]\times[0,d]$ and 
$a$, $b$ and $d$ are real numbers (with $a<b$ and $0<d$).
We assume that the dependency on the parameters is analytic 
and the family is such that the hypotheses {\bf H1} and {\bf H2} 
introduced in section \ref{chapter application}  are satisfied. 
Consider also the reducibility loss bifurcation curves  associated 
to the $2^n$-periodic orbit given by
(\ref{equation alpha n +}). 
Since the value $\alpha_n'(\omega)$ depends also 
on the family of maps $c$ considered, we will denote 
it by $\alpha_n'(\omega,c)$ from now on. We 
omit the case concerning $\beta_n'(\omega)$ (see (\ref{equation alpha n -})) 
because it is completely analogous to the one considered here. 

\subsection{Rotational symmetry reduction}
\label{subsection  rotational symmetry reduction}

Given $\gamma\in \T$,  consider the following auxiliary function
\begin{equation}
\label{definition tgamma}
\begin{array}{rccc}
t_\gamma:&  \BB &\rightarrow &  \BB  \\
\rule{0pt}{5ex} & v(\theta, z) & \mapsto & v(\theta + \gamma, z). 
\end{array}
\end{equation}

Let $\BB_1$ be the subspace of $\BB$ defined by (\ref{equation spaces bbk})
for $k=1$. The space $\BB_1$ is indeed the image of the
projection $\pi_1:\BB \rightarrow \BB$ defined as
\begin{equation}
\label{projection pi_1}
\left[\pi_1 (v)\right](\theta , x)  =  
\left( \int_0^1 v(\theta,x) \cos(2\pi x) d\theta \right) \cos (2\pi\theta) + 
\left( \int_0^1 v(\theta,x) \sin(2\pi x) d\theta \right) \sin (2\pi\theta). 
\end{equation}

Given  $x_0\in \W \cap \R $ and $\theta_0\in\T$ we can
also consider the sets
\[\BB_1' = \BB_1' (\theta_0,x_0) = 
\{f\in \BB_1 \thinspace | f(\theta_0,x_0)=0, \partial_\theta f(\theta_0,x_0)>0 \},\]
and
\[\BB' = \BB' (\theta_0,x_0) = 
\{f\in \BB \thinspace | \pi_1(f) \in \BB_1' \}. 
\]

Note that $ \BB_1' (\theta_0,x_0)$ depends on the election of $(\theta_0,x_0)$, but
for any fixed $x_0$ and $\theta_0$ the set $\BB_1' (\theta_0, x_0)$ is an open
subset of a codimension one linear subspace of $\BB_1$. Note
also that any $v\in$ then $v\neq 0$ due to the condition
$\partial_\theta v(\theta_0,x_0)>0$. Moreover for any $v\in \BB_1 \setminus\{0\}$
there exists  a
unique $\gamma_0\in \T$ such that $t_{\gamma_0}(v)\in \BB_1'(\theta_0,x_0)$
(see proposition 3.1 in \cite{JRT11b}).

Consider a two parametric family of maps $\{c(\alpha,\eps)\}_{(\alpha,\eps)
\in A}$ contained in $\BB$ satisfying the hypotheses {\bf H1} and {\bf H2}
as in section \ref{chapter application}. Consider also the reducibility
loss bifurcation curves  associated to the $2^n$-periodic orbit with
slopes given by (\ref{equation alpha n +}) and (\ref{equation alpha n -}).
Then we have that the formulas
(\ref{equation alpha n +}) and (\ref{equation alpha n -}) can 
be expressed in term of vectors in $\BB_1' (\theta_0, x_0)$. 
Let us see this with more detail. 

Consider the sequences $\{\omega_k\}$, $\{f^{(n)}_k\}$ and $\{u^{(n)}_k\}$ given by
(\ref{equation sequences corollary directions}) and
(\ref{equation sequences directions initial}). Consider now
the sequence
\begin{equation} 
\label{equation v_k in section} 
\tilde{v}^{(n)}_k = t_{\gamma\left(\tilde{v}^{(n)}_{k-1}\right)} 
\left( D \TT_{\omega_{k-1}}  \left(f^{(n)}_{k-1}\right) \tilde{v}^{(n)}_{k-1} \right)
\text{ for } k = 1, ..., n-1, 
\end{equation}
and
\begin{equation}
\label{equation v_0 in section} 
v^{(n)}_0 =  t_{\gamma_0} \left(\partial_\eps c(\alpha_n,0)\right), 
\end{equation}
where $\gamma(\tilde{v}^{(n)}_{k-1})$ and $\gamma_0$ are chosen such
that $\tilde{v}^{(n)}_{k}$ belongs
to $\BB' (\theta_0, x_0)$ for $k=0,1, ..., n$.

If the projection of $D\TT_{\omega_{k-1}}  \left(f^{(n)}_{k-1}
\right) \tilde{v}^{(n)}_{k-1}$ in $\BB_1$ is non zero, then
$\gamma\left(\tilde{v}^{(n)}_{k-1}\right)$ is uniquely determined and
the vectors $\tilde{v}^{(n)}_k$ are well defined. Moreover, 
if we assume that $\omega_0\in \Omega$ and that the 
conjecture {\bf \ref{conjecture H2}} is true then we 
have that $\alpha'_n (\omega)$ can be written as
\begin{equation}
\label{equation alpha n + section}
\alpha'_n (\omega) = 
- \frac{m \left(
DG_1 \left(\omega_{n-1}, f^{(n)}_{n-1}\right) \tilde{v}_{n-1}^{(n)}
\right)
}{\rule{0ex}{3.5ex}  
D \widehat{G}_1 \left(f^{(n)}_{n-1}\right) u_{n-1}^{(n)}} ,
\end{equation}
where $G_1$,  $\widehat{G}_1$ and $m$
are given by equations (\ref{equation definition G1}),
(\ref{equation definition widehatG1}) and (\ref{equation definition of minim}).
For more details see theorem 3.2 in \cite{JRT11b}. 

\subsection{Reduction to the dynamics of the renormalization operator}
\label{subsection reduction to the dynamics}

The goal of this section
is to reduce the problem of describing the asymptotic behavior
of ${\alpha'_n(\omega_0,c_1)}/{\alpha'_{n-1}(\omega_0,c_1)}$ to
the dynamics of the quasi-periodic renormalization operator.

\begin{defin}
\label{definition equivalence of banach sequences}
Given two sequences $\{r_i\}_{i\in \Z_+}$ and $\{s_i\}_{i\in \Z_+}$ in 
a Banach space, we say that they are {\bf asymptotically
equivalent} if there exists $0<\rho <1$ and $k_0$ such that
\[
\| r_i - s_i\| \leq k_0 \rho^i \quad \forall i\in \Z_+ .
\]
We will commit an abuse of notation and
denote this equivalence relation by $s_i \sim r_i$ instead of
$\{r_i\}_{i\in \Z_+} \sim \{s_i\}_{i\in \Z_+}$.
\end{defin}

Given a family $\{c(\alpha,\eps)\}_{(\alpha,\eps)\in A}$ satisfying 
hypotheses  {\bf H1} and {\bf H2} and a fixed Diophantine 
rotation number $\omega_0$, let  $\alpha^*$ denote the parameter value such that the
family $\{c(\alpha,0)\}_{(\alpha,0)\in A}$ intersects
with $W^s(\Phi, \RR)$ and $f_j^*$ denote the intersection
of $W^u(\Phi, \RR)$ with the manifold $\Sigma_j$. Consider then
\begin{equation}
\label{equation sequences directions simplified}
\begin{array}{rcl} 
\rule{0ex}{2.5ex}\omega_k & = & 2 \omega_{k-1}, \text{ for }  k = 1, ..., n-1, \\
\rule{0ex}{4.5ex} 
u_k & = &\displaystyle  \left\{ \begin{array}{ll} 
\rule{0ex}{2.5ex} \displaystyle   D \RR \left(\Phi \right) u_{k-1},  
      &  \text{ for } k = 1, \dots, [n/2]-1, \\ 
\rule{0ex}{2.5ex} \displaystyle   D \RR \left(f_{n-k}^* \right) u_{k-1},  
      &  \text{ for } k = [n/2], \dots, n-1. 
\end{array} \right.  \\
\\\rule{0ex}{4.5ex} 
v_k & = & \displaystyle \left\{ \begin{array}{ll} 
\rule{0ex}{2.5ex} \displaystyle 
t_{\gamma\left(\tilde{v}_{k-1}\right)}
\left( D\TT_{\omega_{k-1}} \left(\Phi \right) v_{k-1} \right),
   & \text{ for } k = 1, ..., [n/2]-1,\\
\rule{0ex}{2.5ex} \displaystyle  t_{\gamma\left(\tilde{v}_{k-1}\right)}
\left( D\TT_{\omega_{k-1}} \left(f_{n-k}^*\right) v_{k-1}\right), 
   &\text{ for } k=[n/2],\dots, n-1.
\end{array} \right. \\
\end{array}
\end{equation}
with
\[  u_0  =  \partial_\alpha c(\alpha^*,0), \quad 
 v_0 =  t_{\gamma_0} \left( \partial_\eps c(\alpha^*,0) \right),   \]
and $\gamma(\tilde{v}_{s-1})$ and $\gamma_0$ are chosen such
that $\tilde{v}^{(n)}_{s}$ belongs
to $\BB_1' (\theta_0, x_0)$ for any $s=1, ..., n$.

Consider the following conjecture for the forthcoming discussion. 

\begin{conj}
\label{conjecture H3}
For any family of maps $\{c(\alpha,\eps)\}_{(\alpha,\eps)\in A}$
satisfying {\bf H1 } and {\bf H2},
assume that
\[
\frac{\tilde{v}^{(n)}_{n-1}}{\|\tilde{v}^{(n)}_{n-1} \|} \sim \frac{v_{n-1}}{\|v_{n-1} \|},
\]
with $\tilde{v}^{(n)}_{n-1}$ and $v_{n-1}$ given by
(\ref{equation v_k in section}) and
(\ref{equation sequences directions simplified}).
 Also assume that there exists a constant $C>0$ such that
\[
\|v_{n-1}\| > C \text{ for any }  n > 0 .  
\]
Finally assume  that there exists a constant $C_0>0$ such that
\[\left|m \left( DG_1\left(\omega_{n-1},f^*_1,
\frac{v_{n-1}}{\|  v_{n-1} \|}\right) \right) \right|> C_0,\]
for any $n\geq0$ and $\omega_0$ Diophantine, where
$m$ is given by (\ref{equation definition of minim}),
$G_1$ by (\ref{equation definition G1}) and $\{f_1^*\}= 
W^u(\RR,\Phi)\cap\Sigma_1 $.
\end{conj}

\begin{table}[t!]
\begin{center}
\texttt{\begin{tabular}{|c|c|c|c|}
\hline 
\rule{0pt}{2.5ex} $n$ & $\|v_{[n/2]-1}\|$ & 
$\left|m \left( DG_1\left(\omega_{[n/2]-1},\Phi,
\frac{v_{[n/2]-1}}{\|  v_{[n/2]-1} \|}\right) \right) \right|$  &   
$\frac{\tilde{v}^{(n)}_{[n/2]-1}}{\left\|\tilde{v}^{(n)}_{[n/2]-1} \right\|} - 
\frac{v_{[n/2]-1}}{\left\|v_{[n/2]-1} \right\|}$ \\
\hline 
4 & 5.426626e+01 & 1.666220e+01 & 9.798573e-02 \\ 
5 & 5.426626e+01 & 1.666220e+01 & 1.019194e-01  \\
6 & 2.361437e+02 & 1.394208e+02 & 4.224092e-03  \\
7 & 2.361437e+02 & 1.394208e+02 & 4.246578e-03  \\
8 & 1.158124e+03 & 5.923824e+02 & 1.410327e-03  \\
9 & 1.158124e+03 & 5.923824e+02 & 1.381662e-03  \\
10 & 6.859354e+03 & 2.901964e+03 & 4.375003e-04  \\ 
11 & 6.859354e+03 & 2.901964e+03 & 4.174722e-04 \\
12 & 5.625187e+04 & 1.717016e+04 & 1.071136e-04 \\
13 & 5.625187e+04 & 1.717016e+04 & 1.024221e-04 \\
\hline
\end{tabular}}
\end{center}
\vspace{-4mm}
\caption{Different values related with the conjecture {\bf \ref{conjecture H3}} 
for the Forced Logistic Map (\ref{FLM new set up}) and $ \omega_0 =  \frac{\sqrt{5}}{2}$. 
}
\label{taula NumSupport conj B}
\end{table}

In the first part of the conjecture we assume that the asymptotic behavior of 
the vectors  $v^{(n)}_{n-1}$ is determined by the linearization of the 
renormalization operator in the fixed point. We have that the iterates of 
$f_0^{(n)}$ correspond to a passage near a saddle point. The 
initial point $f_0^{(n)}$ is always in $\{c(\alpha,0)\}_{(\alpha,0) \in A}$, 
the final point $f_{n-1}^{(n)}$ is 
always in $\Sigma_1$ for any $n$, and the orbit of the points spends more 
and more iterates in a neighborhood of $\Phi$ when $n$ is increased. Therefore 
it is reasonable to expect that the asymptotic behavior is governed by 
the linearization of the operator on the fixed point.  
In the second part of the conjecture we assume that the modulus of the 
vector does not decrease to zero. 

Table \ref{taula NumSupport conj B} supports numerically conjecture 
{\bf \ref{conjecture H3}} for the  Forced Logistic Map (\ref{FLM new set up}) 
with $ \omega =  \frac{\sqrt{5}}{2}$. Note that instead of computing 
the values  $ \|v_{n-1}\| $, $ \left|m \left( DG_1\left(\omega_{n-1},f^*_1,
\frac{v_{n-1}}{\|  v_{n-1} \|}\right) \right) \right|$ and 
$ \frac{\tilde{v}^{(n)}_{n-1}}{\|\tilde{v}^{(n)}_{n-1} \|}  
\frac{v_{n-1}}{\|v_{n-1} \|}$ which appear in conjecture {\bf \ref{conjecture H3}}  
we have computed $\|v_{[n/2]-1}\|$, 
$\left|m \left( DG_1\left(\omega_{[n/2]-1},\Phi,
\frac{v_{[n/2]-1}}{\|  v_{[n/2]-1} \|}\right) \right) \right|$ and 
$\frac{\tilde{v}^{(n)}_{[n/2]-1}}{\left\|\tilde{v}^{(n)}_{[n/2]-1} \right\|} - 
\frac{v_{[n/2]-1}}{\left\|v_{[n/2]-1} \right\|}$. We have done 
this basically to avoid computing the maps $f_i^*= W^u(\Phi,\RR)\cap \Sigma_i$, 
which would require computing $W^u(\Phi,\RR)$ the unstable manifold. 
The point $f_i^*$ accumulate to the fixed point $\Phi$ when
$i\rightarrow \infty$, then one should expect the same behavior 
of the sequences which appear in conjecture {\bf \ref{conjecture H3}}
and the sequence computed in table \ref{taula NumSupport conj B}. For 
more details on how this values are computed see section 
\ref{section applicability FLM}.

Finally we will need the following extension of hypothesis {\bf H2} 

\begin{description}
\item[H3)]   Consider a two parametric family of maps
$\{c(\alpha,\eps)\}_{(\alpha,\eps)\in A}$ (with $A\subset\R^2$)
satisfying
{\bf H1} and {\bf H2} and a fixed Diophantine rotation
number $\omega_0$. Consider also $\omega_n$ and $v_n$ given
by (\ref{equation sequences directions simplified})
and the point  $\{f^*_{1}\} = W^u(\RR,\Phi)\cap \Sigma_1$.
We assume that $DG_1(\omega_{n-1}, f^*_1 )
 v^{(n)}_{n-1}$ has a unique non-degenerate minimum for
any $\omega_0\in \Omega$ and $n\geq 0$. Assume also that
the projection of $D\TT_{\omega_{k-1}}  \left(f^{(n)}_{k-1}
\right) \tilde{v}^{(n)}_{k-1}$ in $\BB_1$ given by
\ref{projection pi_1} is non zero.
\end{description}

Let $\{c(\alpha,\eps)\}_{(\alpha,\eps)\in A}$ (with $A\subset\R^2$) be 
a two parametric family of q.p. forced maps 
satisfying {\bf H1}, {\bf H2} and {\bf H3} and let 
$\omega_0$ be a Diophantine number. Consider the reducibility-loss 
directions $\alpha_n'(c,\omega_0)$ and the sequences $u_n$ and $v_n$ given by 
(\ref{equation sequences directions simplified}). 
In theorem 3.6 of \cite{JRT11b} we proved that, if  
conjectures {\bf \ref{conjecture H2}} and {\bf \ref{conjecture H3}} are true, then
\begin{equation}
\label{equation theorem reduction to dynamics of renormalization}
\frac{\alpha_n'(c,\omega_0)}{ \alpha_{n-1}'(c,\omega_0)} \sim
\deltabf^{-1} \cdot \frac{
\displaystyle m \left( DG_1\left(\omega_{n-1}, f^*_{1} , 
\frac{v_{n-1}}{\|v_{n-1}\|} \right) \right)
}{ \rule{0pt}{4ex} \displaystyle 
m \left( DG_1\left(\omega_{n-2},  f^*_1,  \frac{v_{n-2}}
{\|v_{n-2}\|} \right) \right)
}\cdot 
\left\| D\TT_{\omega_{n-2}}(\Psi) \frac{ v_{n-2}}{\|v_{n-2}\|}
\right\|,
\end{equation} 
where $m$ is given by (\ref{equation definition of minim}),
$G_1$ by (\ref{equation definition G1}), $\{f_1^*\}= 
W^u(\RR,\Phi)\cap\Sigma_1$ is the intersection of the unstable 
manifold of $\RR$ at the fixed point $\Phi$ with the manifold 
$\Sigma_1$ and $\deltabf$ is the universal Feigenbaum constant.

This reduces the asymptotic study of the ratios $\frac{\alpha_n'(c,\omega_0)}
{ \alpha_{n-1}'(c,\omega_0)}$ to the sequence of vectors 
$\frac{ v_{k}}{\|v_{k}\|}$, which are determined by the dynamics of the 
q.p. renormalization operator.

\subsection{Theoretical explanation to the first numerical observation} 
\label{subsection explanation 1st observation}

Consider a two parametric family of maps $\{c(\alpha,\eps)\}_{(\alpha,\eps)
\in A}$ contained in $\BB$ satisfying the hypotheses {\bf H1}, {\bf H2}
and {\bf H3}. Let $\omega_0$ be a  Diophantine rotation number for the family. 

The values $\frac{\alpha'_n(\omega_0,c)}{\alpha'_{n-1}(\omega_0,c)}$ 
depend only on the sequences $\omega_n$ and $v_n$ given by 
(\ref{equation sequences directions simplified}), 
with $v_0=\partial_\eps c(\alpha^*,0) $  and
$\alpha^*$ the parameter value for which the family intersects 
$W^s(\RR,\Phi)$. The behavior of vectors $v_n$ is described by
the dynamics of the following operator,
\begin{equation}
\label{equation map linearized q.p. renormalization general} 
\begin{array}{rccc}
L:& \T \times \BB' &\rightarrow &  \T \times \BB'  \\
\rule{0pt}{5ex} & (\omega, v) & \mapsto & \left(2\omega,  
\displaystyle  \frac{t_{\gamma(v)} \left( D\TT_\omega(\Phi) v\right)}
{\|t_{\gamma(v)} \left(D\TT_\omega(\Phi) v\right) \|} \right),
\end{array}
\end{equation}
where $\gamma$ is chosen such that $t_{\gamma(v)} 
\left( D\TT_\omega(\Phi) v\right)$ belongs to $\BB'$.

In order to study numerically the map $L$ let us recall 
that the spaces $\BB_k$ given by (\ref{equation spaces bbk}) are invariant by 
$D\TT_\omega(\Phi)$, moreover the restriction of $D\TT_\omega(\Phi)$
to this space is equivalent to the map $\LL_{k\omega}$ 
given by (\ref{equation maps L_omega}) with the space $\BB_1$ 
identified with $\RHH(\W) \oplus \RHH(\W)$. For more details 
see proposition 2.16 in \cite{JRT11a}. 

Let us consider the following maps: 
\begin{equation}
\label{equation map linearized q.p. renormalization} 
\begin{array}{rccc}
\tilde{L}_1:& \T \times \BB_1 &\rightarrow &  \T \times \BB_1  \\
\rule{0pt}{5ex} & (\omega, v) & \mapsto & \left(2\omega,  
\displaystyle  \frac{\LL_\omega v}{\|\LL_\omega  v\|} \right) ,
\end{array}
\end{equation}
and 
\begin{equation}
\label{equation map linearized q.p. renormalization section} 
\begin{array}{rccc}
L_1:& \T \times \BB_1' &\rightarrow &  \T \times \BB_1'  \\
\rule{0pt}{5ex} & (\omega, v) & \mapsto & \left(2\omega,  
\displaystyle  \frac{\LL_\omega' v}{\|\LL_\omega'v\|} \right) ,
\end{array}
\end{equation}
with 
\begin{equation}
\label{equation map LL_omega'}
\begin{array}{rccc}
\LL_\omega': &   \BB_1'  & \rightarrow &  \BB_1'   \\
&\rule{0ex}{3ex} v  & \mapsto & \displaystyle t_{\gamma(v)} 
\left( \LL_\omega(v) \right) , 
\end{array}
\end{equation}
where $\gamma(v)$ is chosen such that
$t_{\gamma(v)} \left( \LL_\omega(v)\right) \in \BB_1' $.

We have that $L$ restricted to $\BB_1'$ is equivalent to the map
$L_1$. When $L$ is restricted to $\BB_k$ with $k\neq 1$, then 
it is equivalent to $\tilde{L_1}$. Actually we have that 
the map $L_1$ is equivalent to the map $\tilde{L_1}$ 
after applying the rotational symmetry reduction described in 
section \ref{subsection  rotational symmetry reduction}.

\begin{figure}[t]
\begin{center}
\includegraphics[width=5cm]{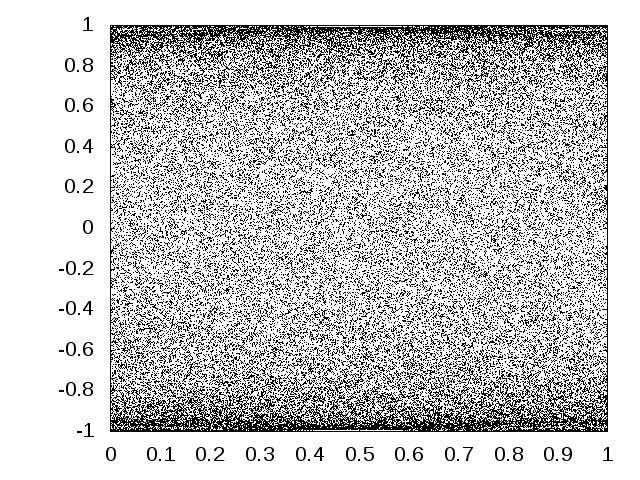}
\includegraphics[width=5cm]{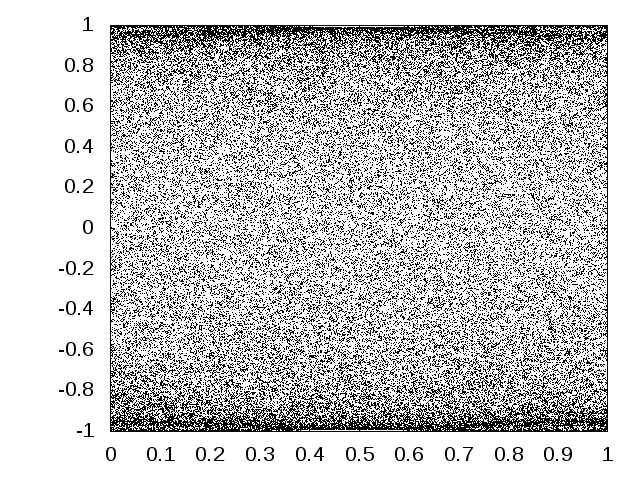}
\includegraphics[width=5cm]{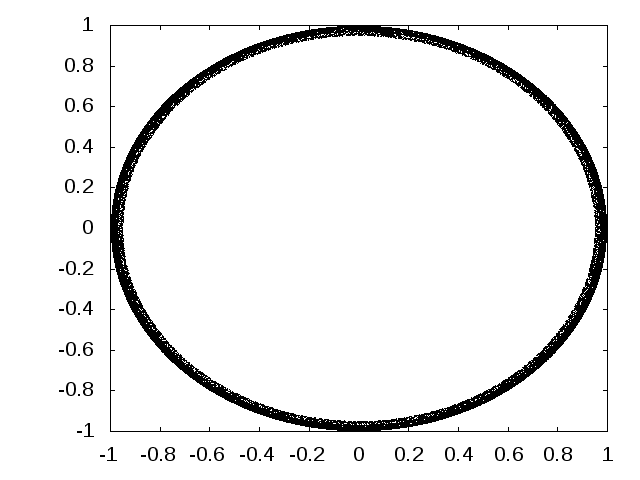}
\includegraphics[width=5cm]{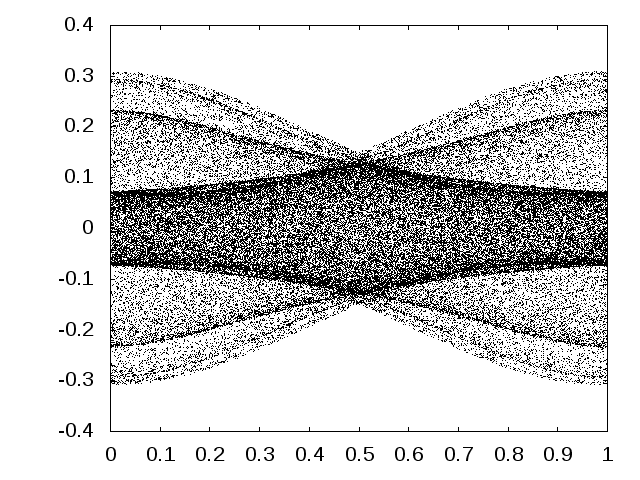}
\includegraphics[width=5cm]{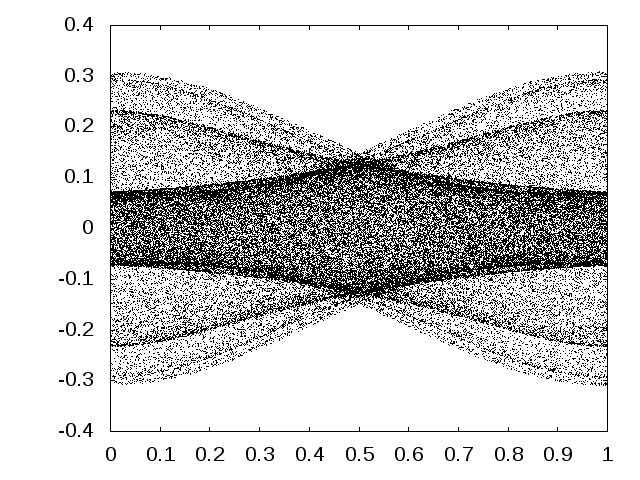}
\includegraphics[width=5cm]{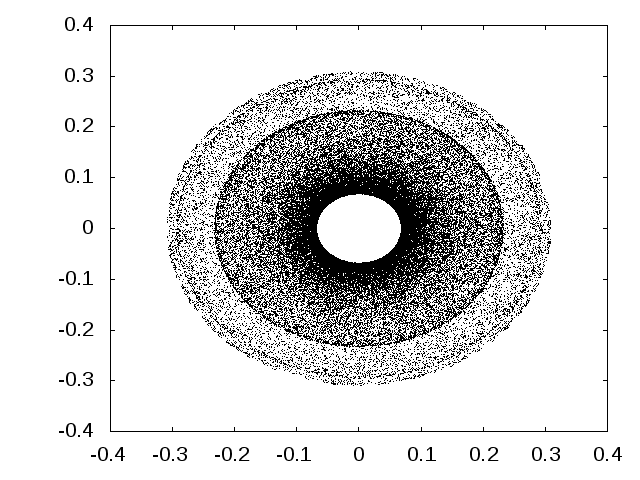}
\includegraphics[width=5cm]{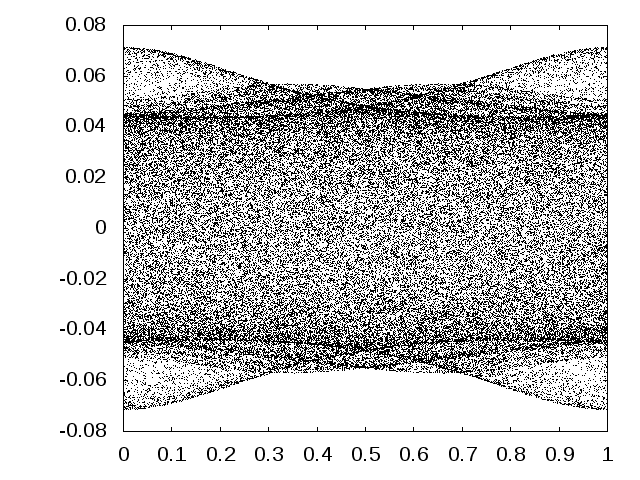}
\includegraphics[width=5cm]{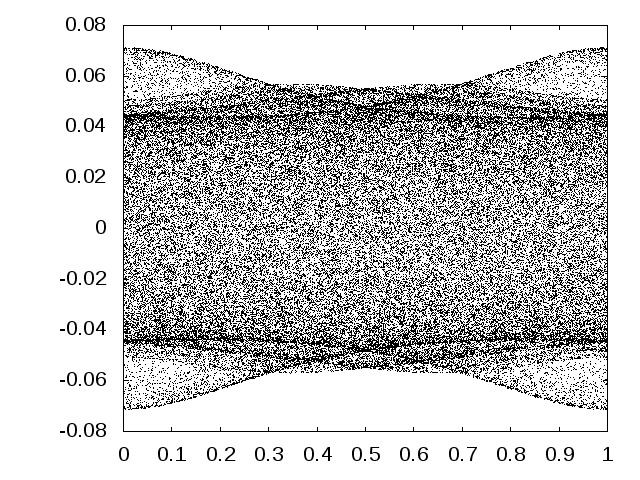}
\includegraphics[width=5cm]{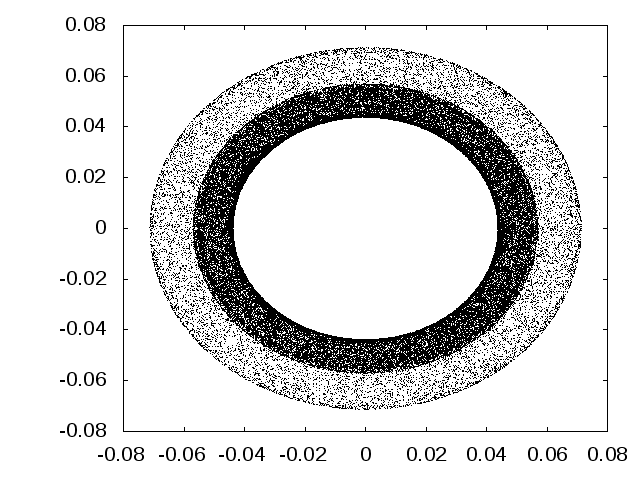}
\end{center}
\caption{
Several planar projections of the attractor of
the map (\ref{equation map linearized q.p. renormalization}).
Form left to right, and top to bottom we have the projections in the
coordinates $(\omega,x_0)$,  $(\omega,y_0)$,  $(x_0,y_0)$,
$(\omega,x_2)$,  $(\omega,y_2)$,  $(x_2,y_2)$,
$(\omega,x_4)$,  $(\omega,y_4)$ and $(x_4,y_4)$. }
\label{figure projections of q.p. renormalization 1}
\end{figure}

\begin{figure}[t!]
\begin{center}
\includegraphics[width=7.5cm]{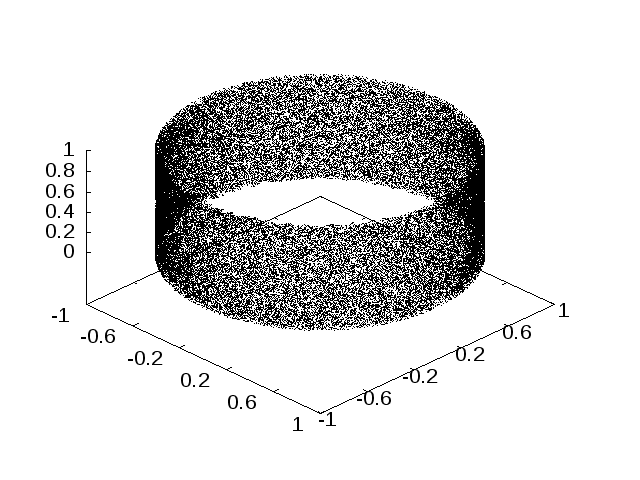}
\includegraphics[width=7.5cm]{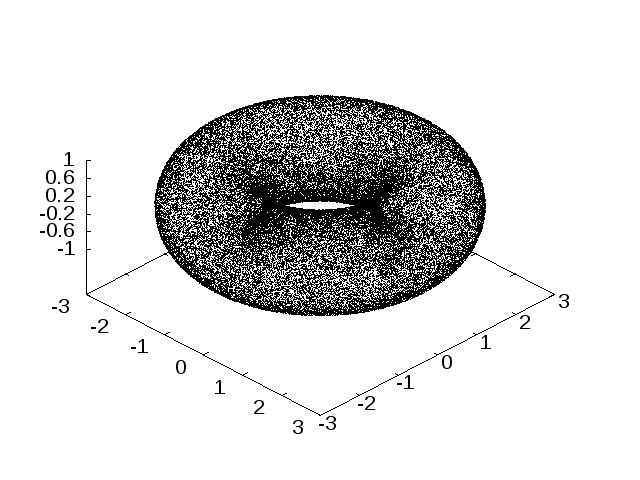}
\includegraphics[width=7.5cm]{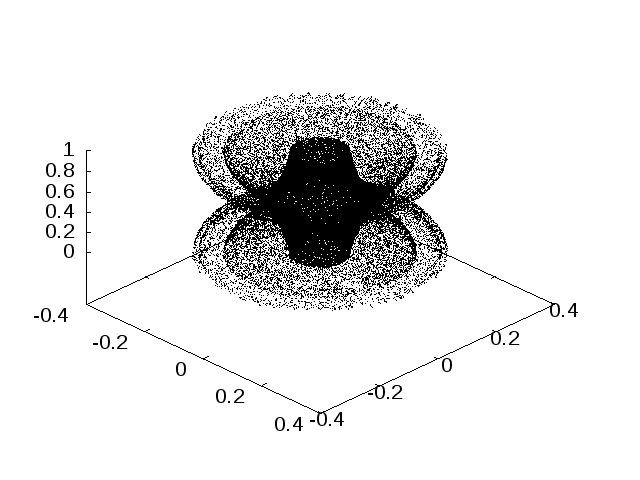}
\includegraphics[width=7.5cm]{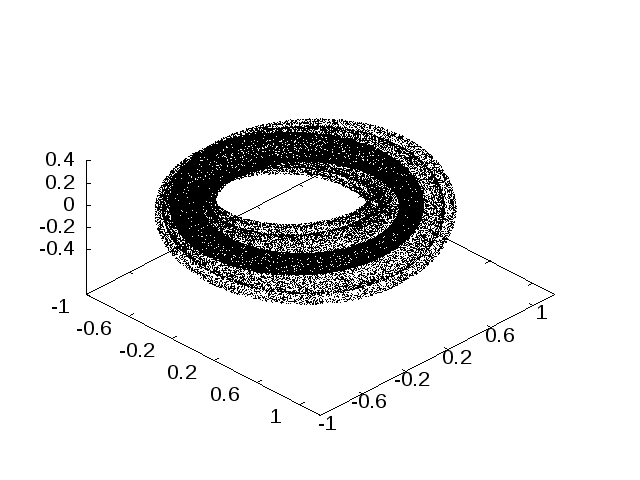}
\includegraphics[width=7.5cm]{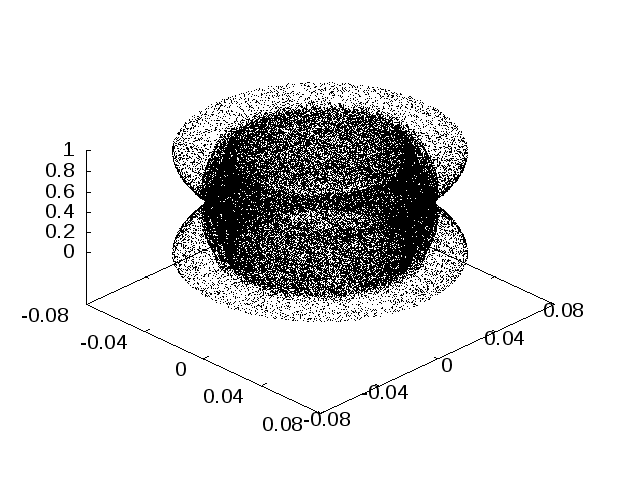}
\includegraphics[width=7.5cm]{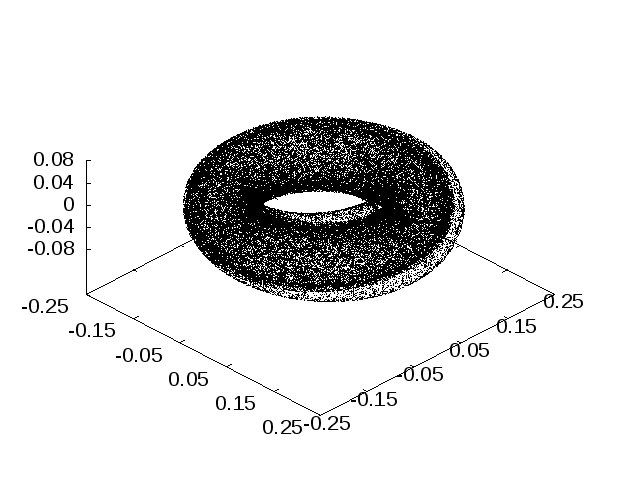}
\end{center}
\caption{
Several spatial projections of the attractor of
the map (\ref{equation map linearized q.p. renormalization}).
In the left hand side of the picture we have a plot (from top to bottom) of the
projections in the coordinates $(\omega,x_0,y_0)$,
$(\omega,x_2,y_2)$ and $(\omega,x_4,y_4)$. In the right hand side there are
displayed the image of the left side projections taking a map that
embeds the solid torus in $\R^3$ (see the text for more details).}
\label{figure projections of q.p. renormalization 2}
\end{figure}

We can use the discretization of $\LL_\omega$ 
described in Appendix 
\ref{section numerical computation of the specturm of L_omega} 
to study numerically the maps $\tilde{L}_1$ and $L_1$. Let us focus 
first on the case concerning $\tilde{L}_1$. Given $v\in \BB_1 = 
\RHH(\W) \oplus \RHH(\W)$, consider the coordinates 
of $v=(x,y)$ given by this splitting. Following the discretization
described in Appendix \ref{section numerical computation of the specturm of L_omega}
each function $x\in \RHH(\W)$ can be approximated by the
vector $(x_0,x_1,x_2,...,x_N)\in \R^{N+1}$ where $x_i$ 
is the $i$-th coefficient of the 
Taylor expansion for $x$ around $0$. This also holds for $y$, 
the second component of $v$. With this procedure each element 
$v$ in $\BB_1$ can be approximated by a vector $(x_0,x_1,...,x_N,y_0,...,y_N)$ 
in $\R^{2(N+1)}$. 

Following the same argument we can approximate $\tilde{L}_1: \T\times \BB_1 
\rightarrow \T \times \BB_1$ by a map 
$\tilde{L}_1^{(N)} : \T \times \R^{2(N+1)} \rightarrow \T \times \R^{2(N+1)}$. 
We can use the discretized map $L_1^{(N)}$ to study the dynamics of $\tilde{L}_1$. 
Given an initial point $v_0=(x^0,y^0)\neq 0$, 
we have iterated the point by the map a certain transient $N_1$ 
and then we have plotted the following $N_2$ iterates. Figures 
\ref{figure projections of q.p. renormalization 1} and 
\ref{figure projections of q.p. renormalization 2} show 
different projections of the attracting set. The values taken for this 
discretization are $N=30$, $N_1=2000$ and $N_2=80000$. We have displayed 
the coordinates corresponding to the first even Taylor 
coefficients of the functions  $x$ and $y$. The odd Taylor 
coefficients obtained were all equal to zero.
This last observation suggest that 
the attractor is contained in the set of even functions (note that 
the subspace of $\BB_1$ consisting  of all the even functions 
is invariant by $\LL_\omega$).

The same computations have been done for bigger values of 
$N$ and  the results are the same. This indicates that the 
set obtained  is stable 
with respect to the order of discretization, therefore it is reasonable
to expect that it is close to the true attracting set of the original system.

Let us remark that we have not made explicit the initial values of
$w_0$ and $v_0$ taken for the computations. Indeed, the
results seem to be independent of these values. We have
repeated this computation taking as initial value of $v_0$
all the elements of the canonical base of the discretized
space $\R^{2(N+1)}$ and we have always obtained the same
results. We have also repeated the computations
for several values of $\omega_0$ obtaining always the same results. 

Given a solid torus $\T \times \D_\rho$, with $\D_\rho$ the
disk of radius $\rho$ in $\C$, we have that the map
$f(\theta,x+i y) = (\cos(\theta) (x+K_0), \sin(\theta)(x+K_0), y)$
embeds this torus in $\R^3$, for any $K_0>\rho$. This embedding
can be used for a better visualization of the spatial projections
of the attracting set of the map (\ref{equation map 
linearized q.p. renormalization section}).  In the right hand
side of figure \ref{figure projections of q.p. renormalization 2}
we have plotted the image by the embedding
of the points on the left hand side. The concrete values
of $K_0$ taken are (from above to below) $2$, $3/4$ and $3/20$.

The numerical approximation of the attractor displayed in figures
 \ref{figure projections of q.p. renormalization 1} and
\ref{figure projections of q.p. renormalization 2}  reveal the 
rotational symmetry of the attractor. This is the rotational symmetry 
described in section \ref{subsection  rotational symmetry reduction}.  

\begin{figure}[t]
\begin{center}
\includegraphics[width=5cm]{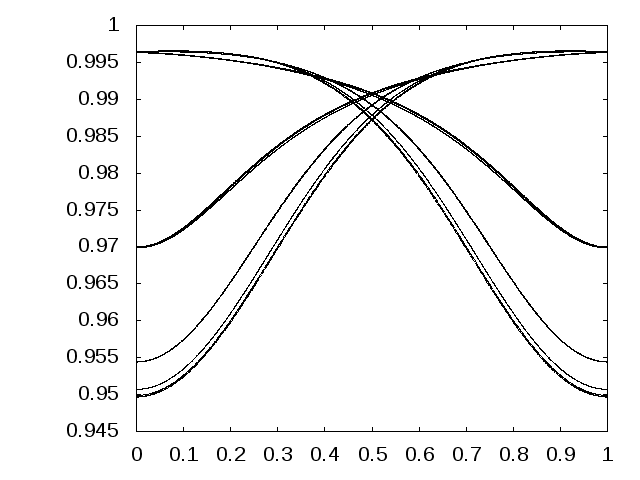}
\end{center}
\begin{center}
\includegraphics[width=5cm]{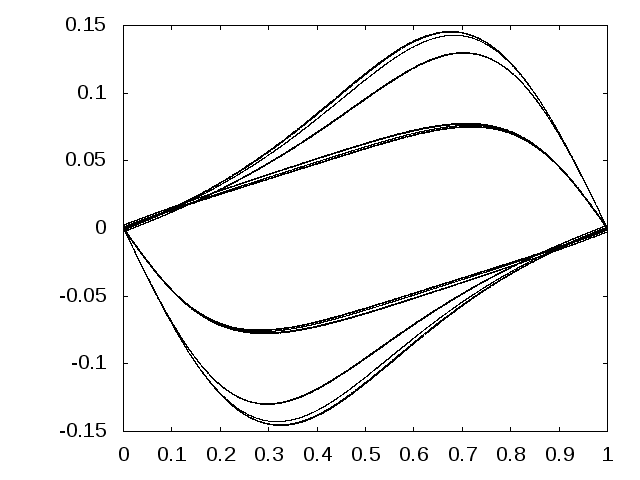}
\includegraphics[width=5cm]{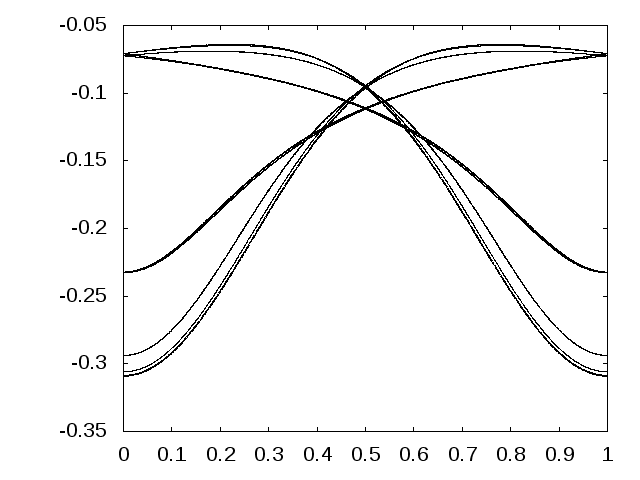}
\includegraphics[width=5cm]{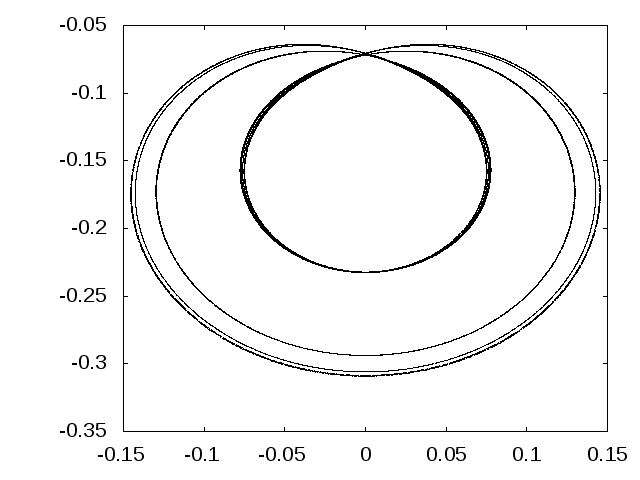}
\includegraphics[width=5cm]{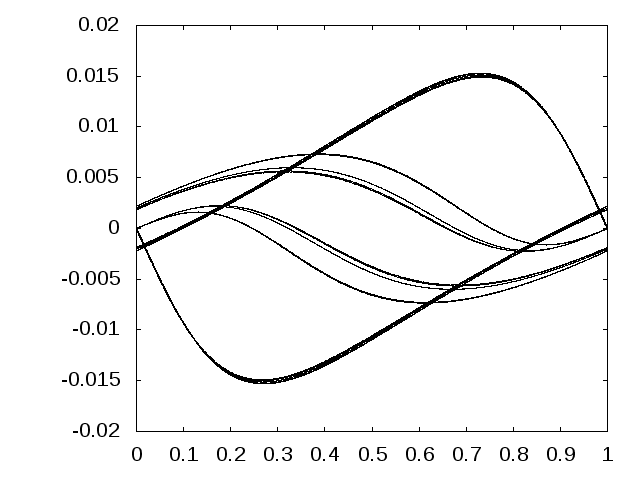}
\includegraphics[width=5cm]{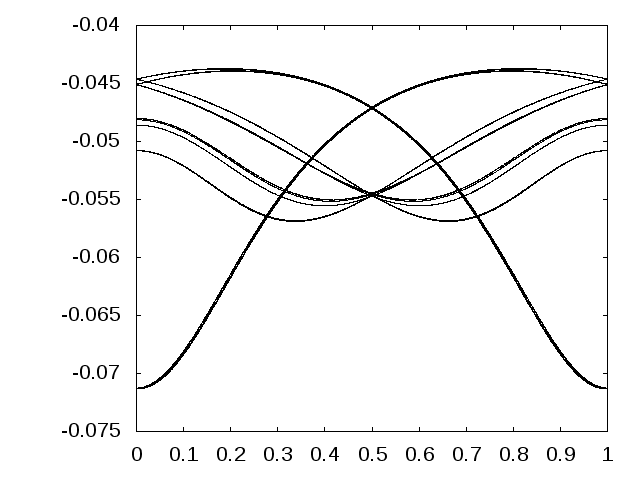}
\includegraphics[width=5cm]{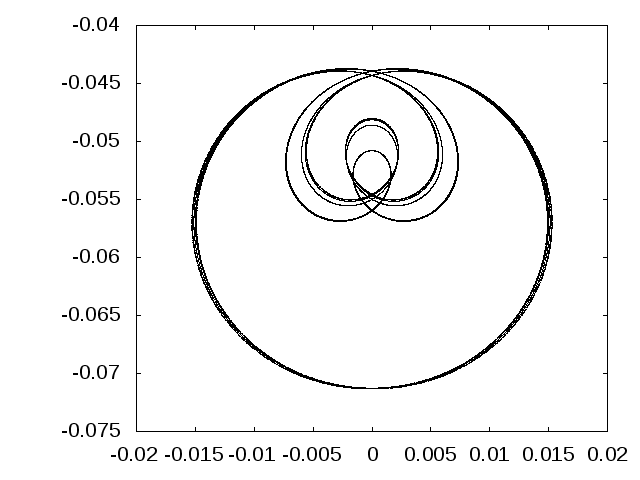}
\end{center}
\caption{
Several planar projections of the section
of attractor of the map (\ref{equation map 
linearized q.p. renormalization section}). Form left to right, and 
top to bottom we have the projections in the
coordinates $(\omega,y_0)$, $(\omega,x_2)$,  $(\omega,y_2)$,  $(x_2,y_2)$,
$(\omega,x_4)$,  $(\omega,y_4)$ and $(x_4,y_4)$. }
\label{figure projections of q.p. renormalization 3}
\end{figure}

\begin{figure}[t]
\begin{center}
\includegraphics[width=7.5cm]{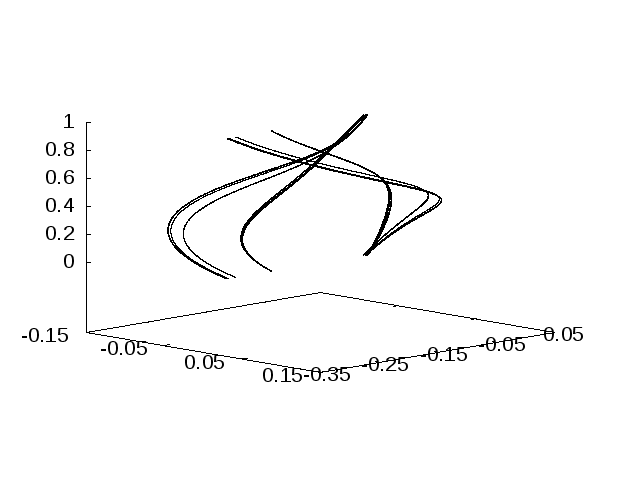}
\includegraphics[width=7.5cm]{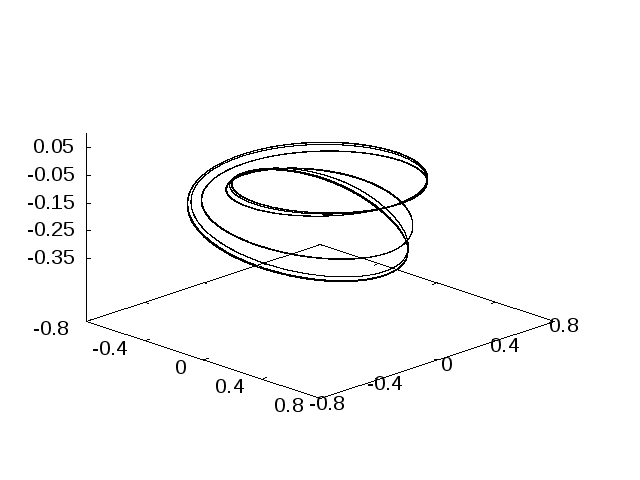}
\includegraphics[width=7.5cm]{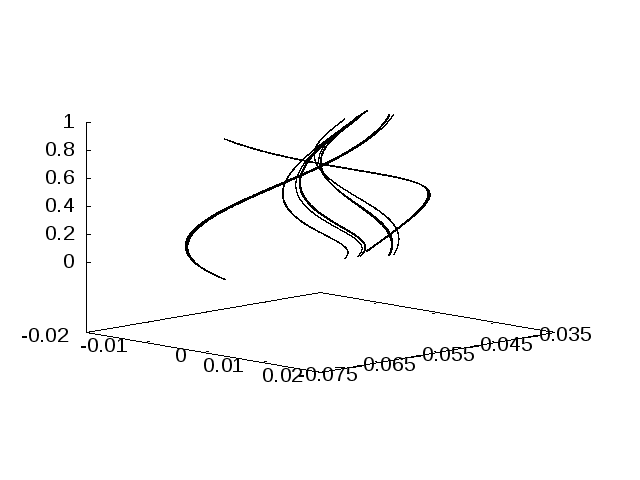}
\includegraphics[width=7.5cm]{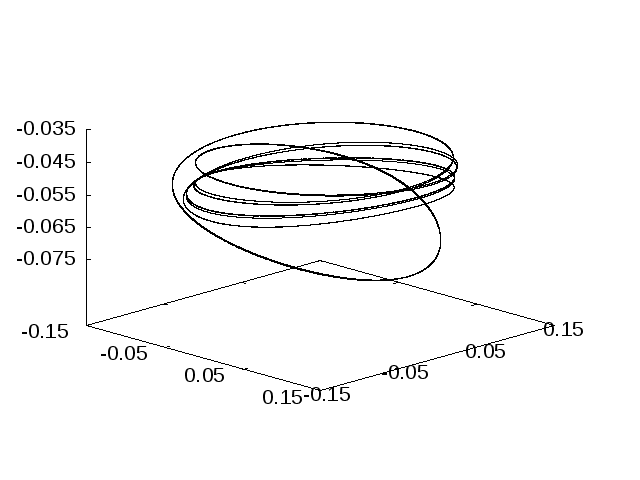}
\end{center}
\caption{
Several spatial projections of the intersection of the attractor of
the map (\ref{equation map linearized q.p. renormalization section}).
The figures in the left  
correspond to the projection to the coordinates $(\omega,x_2,y_2)$ (top)
and $(\omega,x_4,y_4)$ (bottom).
In the right hand side there are
displayed the image of the left side projections taking a map that
embeds the solid torus in $\R^3$ (see the text for more details).}
\label{figure projections of q.p. renormalization 4}
\end{figure}

We can use the discretization described in Appendix 
\ref{section numerical computation of the specturm of L_omega} 
to approximate numerically the dynamics of $L_1$, 
as we have just done  in the case of $\tilde{L}_1$. For the numerical simulation 
of the operator, we have taken $\theta_0=0$ and $x_0=0$. In this case 
the set $\BB'_1(0,0)$ is identified in $\R^{2(N+1)}$ with the 
half hyperplane 
\[
\{(x,y)\in \R^{2(N+1)} \thinspace | \thinspace x_0=0 \text{ and } y_0> 0\},
\]
where $x_0$ and $y_0$ are respectively  the first components of $x$ and $y$. 

In figures \ref{figure projections of q.p. renormalization 3} and
\ref{figure projections of q.p. renormalization 4} there are displayed
different projections of the attracting set obtained 
iterating the map $L_1$.  
As before we have also considered the map $f(\theta,x+i y) =
(\cos(\theta) (x+K_0), \sin(\theta)(x+K_0), y)$ which embeds 
the solid torus $\T \times \D_\rho$ in $\R^3$ for a better
visualization of the set. This
time the values of $K_0$ have been taken equal to
$1/2$ (above) and $3/25$ (below).

Note that the different projections of the attracting set 
displayed in figure \ref{figure projections of q.p. renormalization 3} 
keep a big resemblance with the plots of the dyadic solenoid displayed 
in figure 5 of \cite{Mil85}. 
Indeed we believe that the attractor is the inclusion of a  dyadic 
solenoid in $\BB_1'$. For more details on the definition 
and the dynamics of the solenoid map see \cite{BT11, KH96, Mil85, Sma67}. 
To explain this fact, let us introduce a new 
conjecture, this time on the operator $L_1$.

\begin{figure}[t]
\begin{center}
\includegraphics[width=5cm]{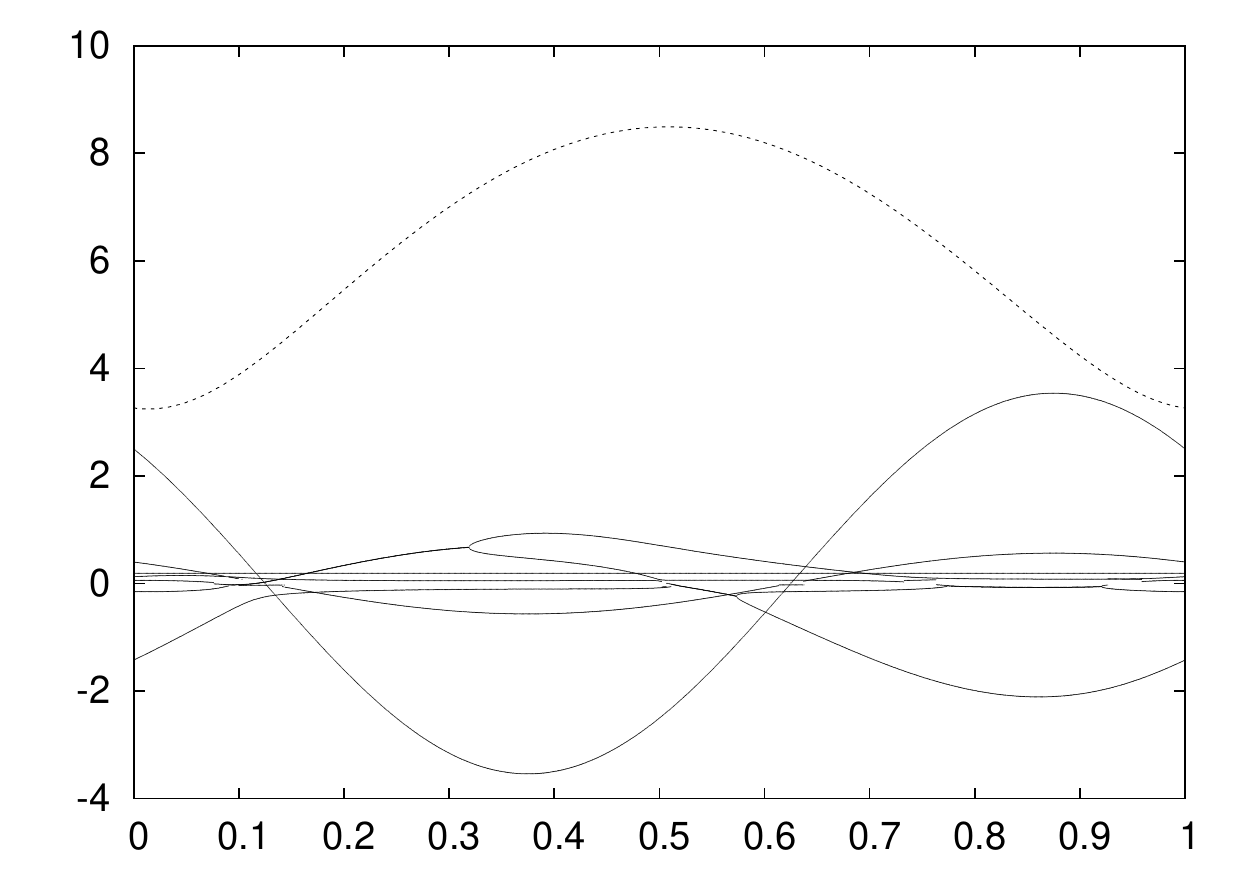}
\includegraphics[width=5cm]{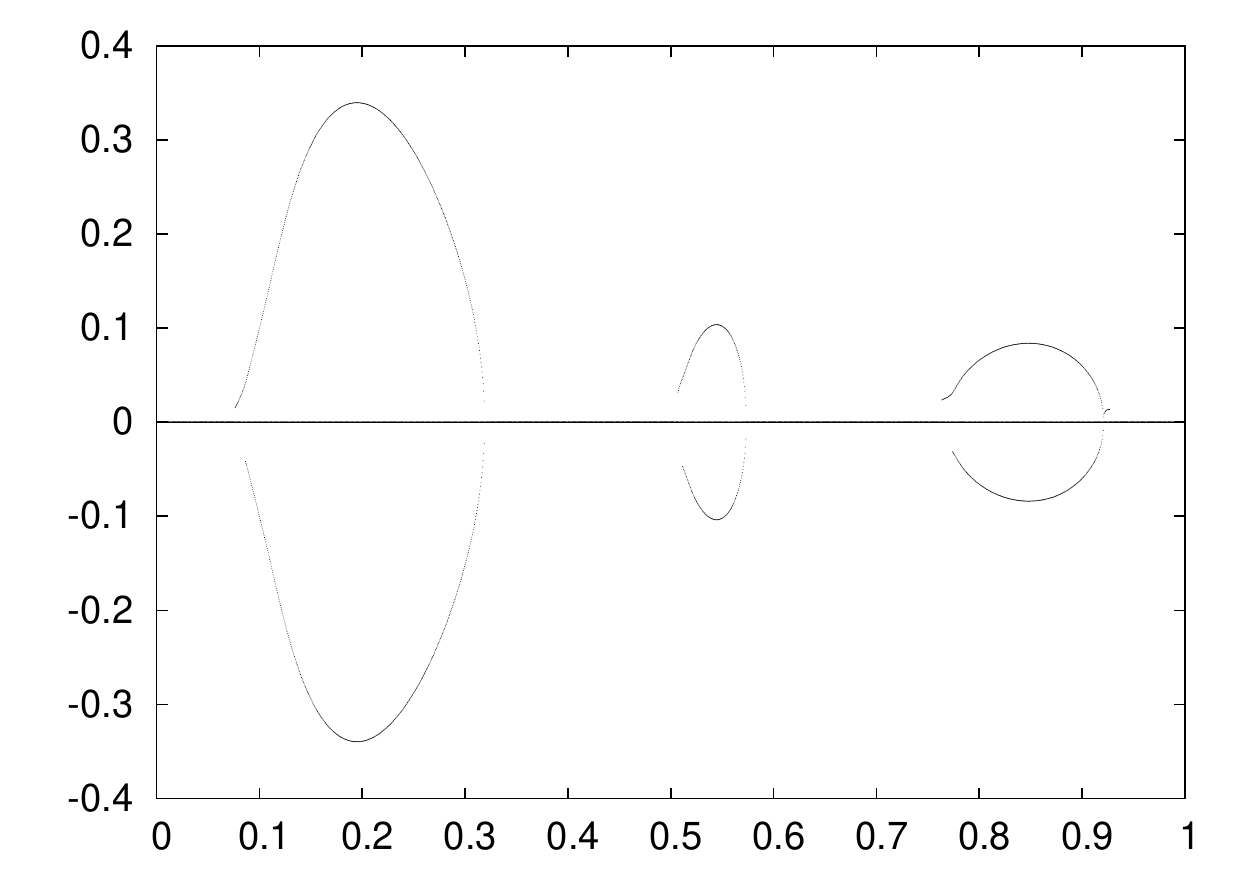}
\includegraphics[width=5cm]{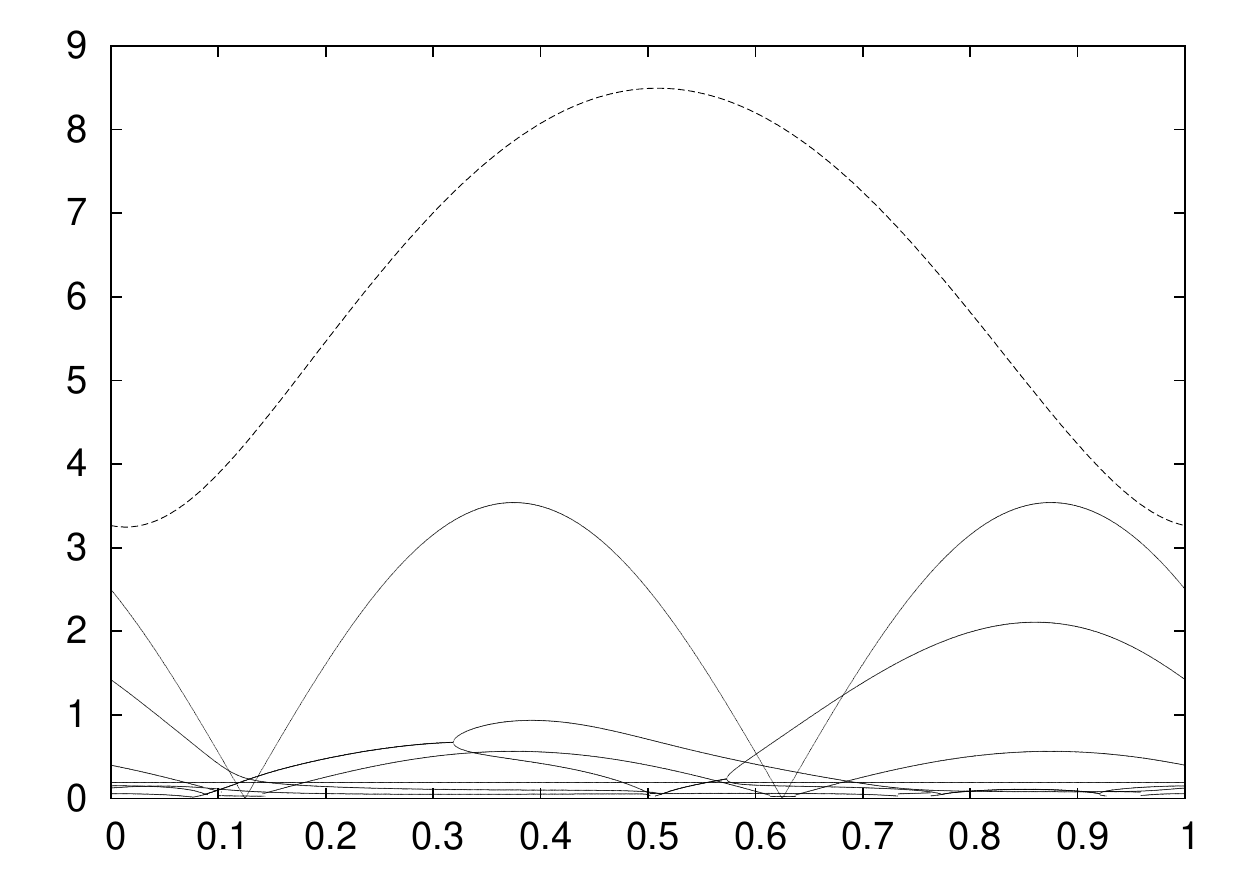}
\end{center}
\caption{
Numerical approximation of the spectrum of 
$\LL_\omega'$ with respect to the parameter $\omega$. 
From left to right we have the real part, 
the imaginary part and the modulus of the first eight eigenvalues
 of $\LL_\omega'$ with respect to $\omega$. }
\label{figure espectrum of L_omega'}
\end{figure}

\begin{conj}
\label{conjecture H4}
There exist an open set $V\subset 
\BB_1'$ (independent of $\omega$)  such that 
the second component of the map $L_1$ given by 
(\ref{equation map linearized q.p. renormalization section}) is 
contractive (with the supremum norm) 
in the unit sphere and it maps the set $V$ into itself 
for any $\omega \in \T$. Additionally we will assume that the contraction 
is uniform for any $\omega \in \T$, in the sense that there exists 
a constant $0<\rho<1$ such that the Lipschitz constant 
associated to the second component 
of the map $L_1'$ is upper bounded by $\rho$ for any $\omega\in \T$. 
\end{conj}

A good reason to think that conjecture {\bf \ref{conjecture H4}} is
true resides in the spectrum of the operator $\LL'_\omega$.
In figure \ref{figure espectrum of L_omega'} we have
a numerical approximation of this spectrum with respect to
$\omega$ as a parameter.  For this
computation we have followed the same procedure that
we used for the computation of the
spectrum of $\LL_\omega$, for details see
Appendix \ref{section numerical computation of the specturm of L_omega}.

Looking at figure \ref{figure espectrum of L_omega'} we can observe 
that there exists a dominant eigenvalue (which is plotted in a dashed 
line) that does not cross the rest, which varies ``nicely'' with 
respect to $\omega$. Then for each value of $\omega$, 
the normalization of $\LL'_\omega$  is a contraction in the sphere, 
with the eigenvector associated to the dominant eigenvalue as a fixed point. 
This means that conjecture {\bf \ref{conjecture H4}}  is true ``point-wise'', 
but this is not enough because the domain of contractivity might depend
 on $\omega$. 

Let us justify now why conjecture {\bf \ref{conjecture H4}} would explain 
the numerical results obtained for the attractor of $L_1$.  Consider the set
$\T\times V \subset \T \times \BB'_1(x_0,\theta_0)$. If the conjecture
{\bf \ref{conjecture H4}}
is true, then we would have
that the set $v$ invariant by the map $L_1$. More concretely
we would have
that the set would be expanded to the double of its length on
the $\T$ direction and contracted in the $\BB_1'$ direction.
Assume that this transformation is done in such a way that
$L_1$ maps the set $\T\times V$ inside itself but without
self intersections. When we consider the
intersection of $L_1(\T\times V)$ with a set of the type
$\{\gamma_0\}\times V$ (for some $\gamma_0\in T$) the
section is conformed by two different sets without
self intersections. The subsequent images by $L_1$ we would have
(for each leaf $\{\gamma_0\}\times V$) the double of
components than in the previous step, each
of them strictly contracted in the $\BB'_1$ component
and contained in the previous set. Note that the described process is 
completely analogous to the geometric construction of a dyadic solenoid,
but this time contained inside the Banach space $\BB'_1$ instead 
of the solid torus. Therefore conjecture {\bf \ref{conjecture H4}} would give 
an explanation to the numerical approximation of the attracting set of 
$L_1$ obtained before.

Consider $\{c_1(\alpha,\eps)\}$ and $\{c_2(\beta,\eps)\}$
two different families of two parametric
maps satisfying hypotheses {\bf H1}, {\bf H2'} and {\bf H3}.
Since the  family of maps satisfy hypothesis  {\bf H2'}, 
we have that $\partial_\eps c_i(\alpha^*,0)$  belongs to $\BB_1$ for $i=1,2$. 
Therefore the dynamics of $L$ 
(\ref{equation map linearized q.p. renormalization general})  coincide with 
the dynamics of $L_1$ (\ref{equation map linearized q.p. renormalization section}).
Theorem 3.10 in \cite{JRT11b} asserts the following. 

Consider $\{c_1(\alpha,\eps)\}$ and $\{c_2(\beta,\eps)\}$
two different families of two parametric
maps satisfying hypotheses {\bf H1} and {\bf H2'}.
Let $\alpha^*$ and $\beta^*$ be the parameter
values where each family $c_1(\alpha,0)$ and $c_2(\beta,0)$
intersects $W^s(\RR,\Phi)$, the
stable manifold of the fixed point of the renormalization operator.
Let $ \operatorname{Rot}(V) = \left\{ v \in \BB_1 \thinspace
 | \thinspace t_\gamma(v) \in V \subset \BB_1' 
\text{ for some } \gamma \in \T\right\} $ where $V$ is the set
given by conjecture {\bf \ref{conjecture H4}} for some $\gamma\in\T$.
If $\partial_\eps c_i(\alpha^*,0)$  belongs to $\operatorname{Rot}(V)$ 
for $i=1,2$ and the conjectures {\bf \ref{conjecture H2}}, 
{\bf \ref{conjecture H3}} and {\bf \ref{conjecture H4}} are true, 
then (for any $\omega_0\in \Omega$)  we have that
\begin{equation}
\label{equation thm proof universality conjecture} 
\frac{\alpha_n'(\omega_0,c_1)}{ \alpha_{n-1}'(\omega_0,c_1)}
\sim
\frac{\alpha_n'(\omega_0,c_2)}{ \alpha_{n-1}'(\omega_0,c_2)},
\end{equation}
where $\alpha_i'(\omega_0, c_i)$  are the reducibility loss directions 
associated to each family $c_i$ for the rotation number of the system 
equal to $\omega_0$.

\subsection{Theoretical explanation to the second numerical observation} 
\label{subsection explanation 2nd observation}

As before consider a two parametric family of maps $\{c(\alpha,\eps)\}_{(\alpha,\eps)
\in A}$ satisfying hypothesis {\bf H1}, {\bf H2'} and {\bf H3'} 
and let 
$\alpha'_n(\omega, c)$ denote the slope of one of the curves 
of the reducibility loss bifurcation 
associated to the $2^n$ periodic invariant curve of the family. 
In theorem 3.13 of \cite{JRT11b} we prove that the second numerical 
observation done in section \ref{subsection intro numerical evidences} can be
explained as a consequence of  the universal behavior 
(\ref{equation thm proof universality conjecture}). One 
of the hypothesis to prove this requires $\alpha'_i(\omega,c)/
\alpha'_i(2\omega,c)$ to be a bounded sequence.

The boundedness of this sequence can be obtained on its turn as 
a consequence of the following conjecture on the operator $\LL_\omega$.

\begin{conj}
Consider  $\LL_\omega: \RHH(\W_{\rho})
\oplus \RHH(\W_{\rho})\rightarrow \RHH(\W_{\rho}) \oplus \RHH(\W_{\rho})$
the map given by (\ref{equation maps L_omega}) and $\omega_0 \in \Omega$. 
Given  $v_{0,1}$ and $v_{0,2}$ two vectors in $\RHH(\W_{\rho})
\oplus \RHH(\W_{\rho})\setminus\{0\}$, consider the sequences 

\begin{equation}
\label{equation sequences directions counterexample 0}
\begin{array}{rcll} 
\omega_k & = & 2 \omega_{k-1} 
&
\text{ for } k = 1, ..., n-1.\\
\rule{0ex}{4ex} 
v_{k,1} & = & \LL_{\omega_{k-1}} \left(u_{k-1,1} \right)
&  \text{ for } k = 1, ..., n-1.\\
\rule{0ex}{4ex} 
v_{k,2} & = & \LL_{2\omega_{k-1}}\left(u_{k-1,2} \right)
&  \text{ for } k = 1, ..., n-1.
\end{array}
\end{equation}

\label{conjecture H5}
Then, there exist constants $C_1$ and $C_2$ such that
\[
C_1 \frac{\|v_{0,2} \|}{\|v_{0,1}\|} \leq \frac{\|v_{n,2} \|} {\|v_{n,1}\|}  
\leq C_2 \frac{\|v_{0,2} \|}{\|v_{0,1}\|}
\]
for any $n\geq0$.
\end{conj}

\begin{figure}[t]
\begin{center}
\includegraphics[width=7.5cm]{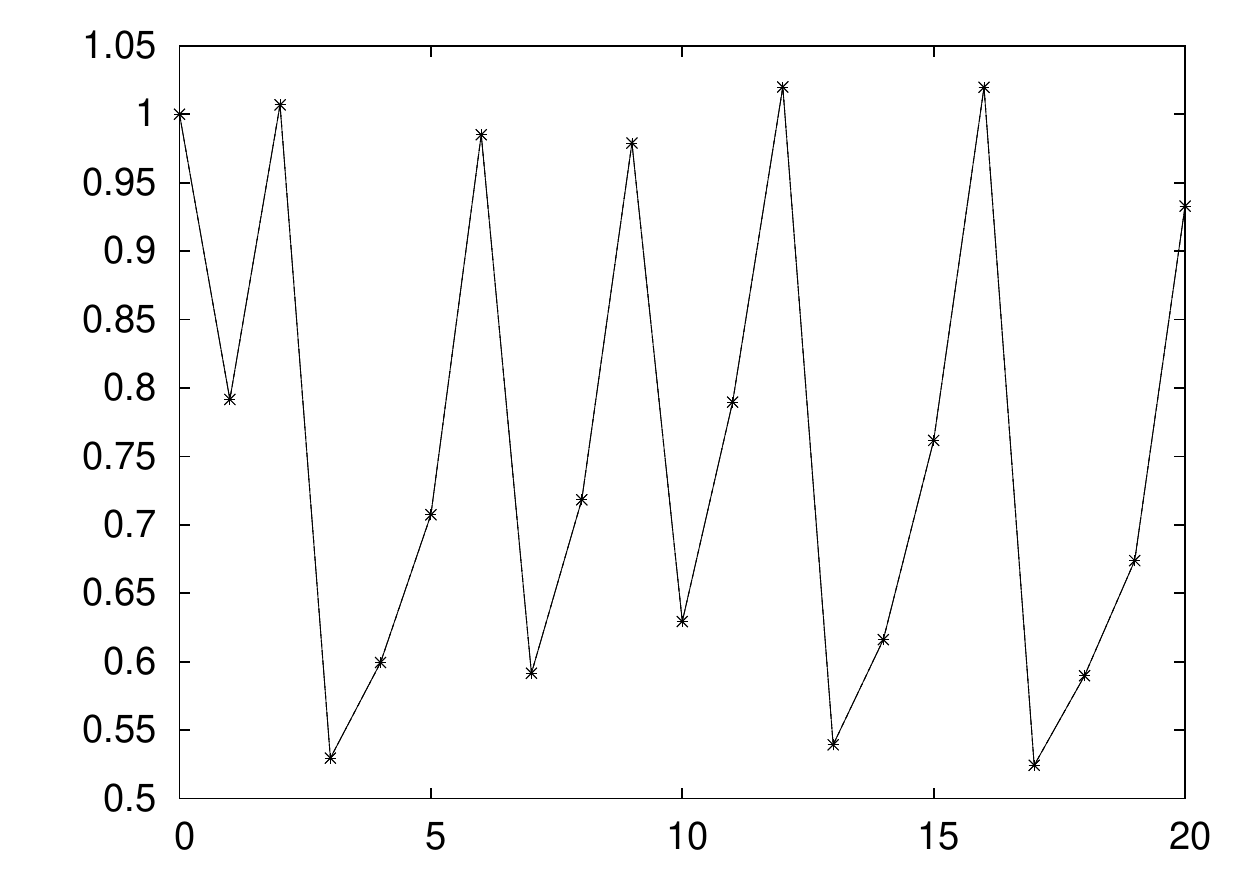}
\includegraphics[width=7.5cm]{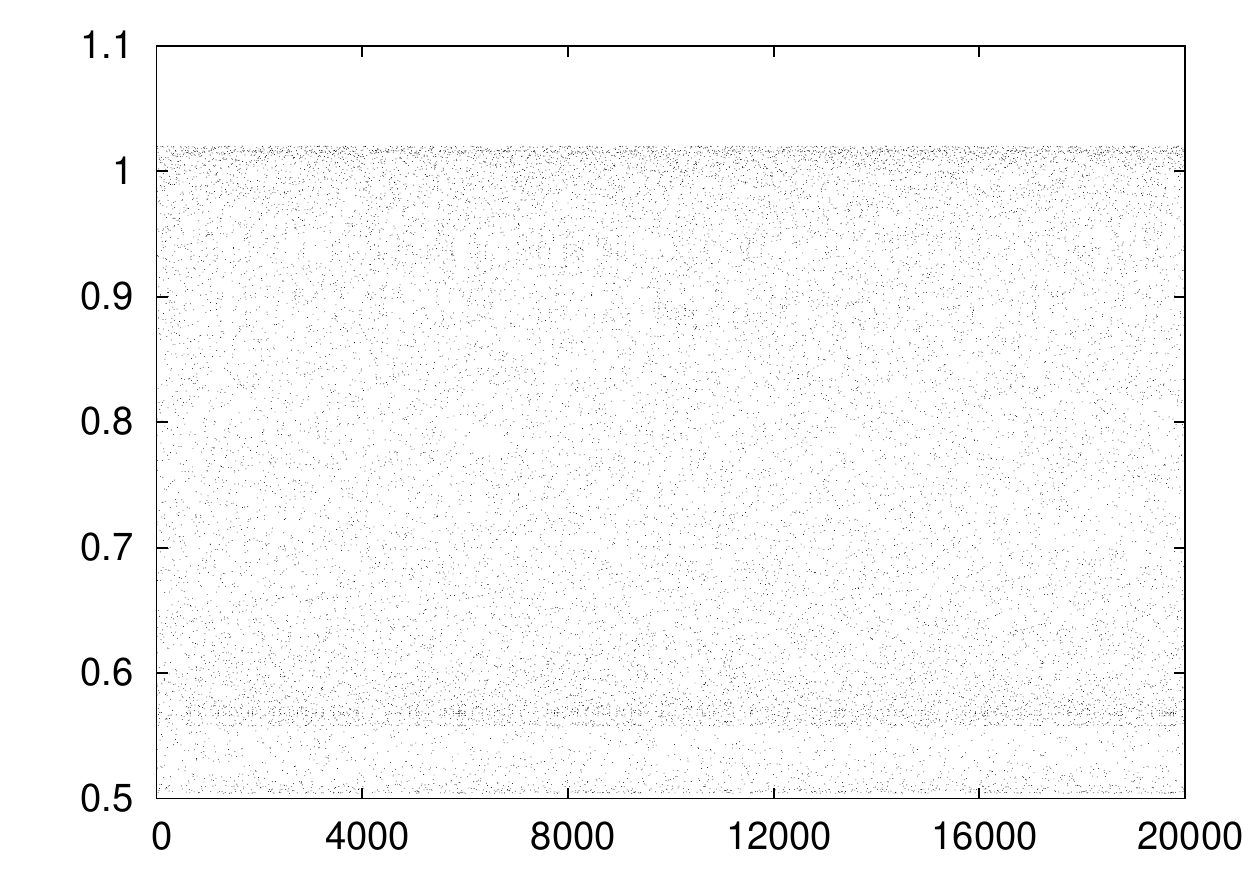}
\end{center}
\caption{ 
We have $\frac{\|v_{0,1}\|}{\|v_{0,2} \|} \frac{\|v_{n,2} \|} {\|v_{n,1}\|}$ 
with respect to $n$, from  $n=0$ to $n=20$ in the left and from
$n=0$ to $n=2 \times 10^4$ in the right.}
\label{figure normes contraexemple}
\end{figure}

This conjecture can be interpreted as a uniform growth condition 
on $\LL_\omega$. To support this conjecture we have computed 
numerically the iterates (\ref{equation sequences directions counterexample 0}), 
in order to estimate the growth of $\|v_{n,1} \|$ with respect to the
growth of $\|v_{n,2}\|$.  
In figure \ref{figure normes contraexemple} we have plotted the ratios
$\frac{\|v_{0,1}\|}{\|v_{0,2} \|} \frac{\|v_{n,2} \|} {\|v_{n,1}\|}$ with respect
to $n$ for the sequence (\ref{equation sequences directions counterexample 0}) with 
$v_{0,1}= 1 $ and $v_{0,2} = 1$. For other initial vectors we obtain the same 
behavior. This suggests that conjecture {\bf \ref{conjecture H5}} is true.

Given $\{c(\alpha,\eps)\}$ a two parametric maps satisfying hypotheses 
{\bf H1}, {\bf H2'} and {\bf H3} let  $\alpha^*$ be the parameter
values for which the family $\{c(\alpha,0)\}_{(\alpha,0)\in A}$
intersects $W^s(\RR,\Phi)$.
Let $ \operatorname{Rot}(V) = \left\{ v \in \BB_1 \thinspace
 | \thinspace t_\gamma(v) \in V \subset \BB_1' 
\text{ for some } \gamma \in \T\right\} $ where $V$ is the set
given by conjecture {\bf \ref{conjecture H4}} for some $\gamma\in\T$.
Assume that $\partial_\eps c(\alpha^*,0) \in \operatorname{Rot}(V)$.
In corollary 3.14 of \cite{JRT11b} we prove that, for any $\omega_0\in \Omega$, 
then 
\begin{equation}
\label{equation coro univ implies renor} 
\lim_{n\rightarrow \infty} 
\frac{\alpha_n'(\omega_0,c)}{ \alpha_{n-1}'(2\omega_0,c)}
\end{equation}
exists, which explains the second numerical observation of
section \ref{subsection intro numerical evidences}.

\subsection{Theoretical explanation to the third numerical observation} 
\label{subsection explanation 3rd observation}

In sections \ref{subsection explanation 1st observation} and 
\ref{subsection explanation 2nd observation} we focussed the
discussion on the asymptotic behavior for families satisfying hypothesis
{\bf H2'}. The map considered in the third numerical observation of 
section \ref{subsection intro numerical evidences} is an example of a
map not satisfying {\bf H2'} for which equations 
(\ref{equation thm proof universality conjecture}) and 
(\ref{equation coro univ implies renor}) do not hold. 

Let $\{c(\alpha,\eps)\}_{(\alpha,\eps) \in A}$  be a two parametric family of maps 
satisfying hypothesis {\bf H1}, {\bf H2} and {\bf H3}. Let
$\alpha'_n(\omega, c)$ denote slope of the reducibility loss bifurcation
associated to the $2^n$ periodic invariant curve of the family. 
Finally consider $\omega_0$ a  Diophantine rotation
number for the family. Let $\alpha^*$ be the parameter value 
for which $\{c(\alpha,0)\}_{(\alpha,0)\in A}$ intersects $W^s(\Phi,\RR)$.

The main difference between the example considered in the third numerical observation of
section \ref{subsection intro numerical evidences} and 
the previous two is that  
\[\partial_\eps c(\alpha^*,0) = v_{0,1} + v_{0,2} \text{ with } 
v_{0,i} \in \BB_i, \quad i=1,2,\] where
the spaces $\BB_i$ are given by (\ref{equation spaces bbk}). 
More concretely, for the numerical example cited above we have 
considered 
\begin{equation}
\label{eq: initial vectors}
v_{0,1}= f_1(x) \cos(\theta), \quad v_{0,2} = \eta f_2(x) \cos(2\theta). 
\end{equation}
As this family depends on $\eta$, we denote by $c_\eta$ this 
concrete family. The parameter $\eta$ is considered in addition 
to the parameters $\alpha$ and $\eps$ of the family. In other 
words, for each $\eta\geq 0$, $c_\eta$ is a two parametric family of maps. 

Then one has that equation 
(\ref{equation theorem reduction to dynamics of renormalization}) still holds 
with $v_n= v_{n,1} + v_{n,2}$, with  $v_{n,1}$ and $v_{n,2}$ the 
sequence given by (\ref{equation sequences directions counterexample 0}) 
and  $v_{0,1}$ and $v_{0,2}$ given by (\ref{eq: initial vectors}). In
this case there is no universal behavior because the sequence $v_n$ 
lives in a bigger invariant space, where the renormalization 
operator is not contractive. On the other hand, the 
numerical computations in \cite{JRT11p2} suggest that the sequence 
$\alpha_n'(\omega_0,c_\eta)/ \alpha_{n-1}'(\omega_0,c_\eta)$ (for $\eta>0$)
is not asymptotically equivalent to $\alpha_n'(\omega_0,c_0)/ 
\alpha_{n-1}'(\omega_0,c_0)$, but both sequences are $\eta$-close to 
each other. This can be explained as a consequence conjecture 
{\bf \ref{conjecture H5}} where we conjecture uniform growth (in norm) of the
sequences $v_{n,1}$ and $v_{n,2}$. For more details see 
theorem 3.20 in \cite{JRT11b}. 
 
\section{Applicability to the Forced Logistic Map} 
\label{section applicability FLM}

The theory exposed in sections \ref{section review q-p renor} 
and \ref{section A not yet rigorous explanation} have 
been built as a response to the observations done in the 
study of the Forced Logistic Map (see \cite{JRT11p, JRT11p2}). 
In this section we discuss the applicability of the 
quasi-periodic renormalization theory to the Forced 
Logistic Map.  In the cited papers we considered two different version of 
the FLM, which correspond to the map (\ref{q.p- family}) with either 
$h_{\alpha,\eps}(\theta,x) = \alpha x (1-x) \cos(2\pi \theta)$  or 
$h_{\alpha,\eps}(\theta,x) = \cos (2 \pi \theta)$, where 
the parameters $(\alpha, \eps) \in [0,4] \times [0,1]$. 
Notice that these two forms of the FLM do not satisfy
the requirements of the theory developed in the
previous sections because the family of maps does not belong to $\BB$.

This problem can be easily solved as follows. For $\alpha >2$ we can 
consider the affine change of variables given by
$y = a x + b$, with $a= \frac{4}{\alpha-2}$ and $b= - \frac{2}{\alpha-2}$. If
we apply this change of variables to the family (\ref{q.p- family}) when 
$h_{\alpha,\eps}(\theta,x) = \alpha x (1-x) \cos(2\pi \theta)$ we obtain 
the following family,
\begin{equation}
\label{FLM new set up}
\left.
\begin{array}{rcl}
\bar{\theta} & =  & \theta + \omega  ,\\
\bar{y} & = & \alpha \displaystyle \left(  
\frac{\alpha}{\alpha -2} - \frac{\alpha(\alpha-2)}{4} y^2 
\right)  (1+ \eps \cos(2\pi \theta)) - \frac{2}{\alpha -2}  . 
\end{array}
\right\}
\end{equation}

If we apply the same change of variables when $h_{\alpha,\eps}(\theta,x) 
= \cos(2\pi \theta)$ we obtain this other family 
\begin{equation}
\label{FLM B new set up}
\left.
\begin{array}{rcl}
\bar{\theta} & =  & \theta + \omega  ,\\
\bar{y} & = & \displaystyle 
 1 - \frac{\alpha(\alpha-2)}{4} y^2   + \frac{4\eps}{\alpha -2} \cos(2\pi \theta). 
\end{array}
\right\}
\end{equation}

Although the change of variables considered depends on the parameter $\alpha$,
the parameter space of the maps (\ref{FLM new set up}) and 
(\ref{FLM B new set up}) is the same as the parameter space of
the map (\ref{q.p- family}) (for the corresponding value of $h$). 
Then any conclusion drawn on the parameter 
space of the map  (\ref{FLM new set up}) (respectively
(\ref{FLM B new set up})) extends automatically to the parameter space of
(\ref{q.p- family}).

With this new set up, we have that both families of maps belong to 
$\BB$ for $\alpha \in (2,4)$ and $\eps$ small enough.
One should check that the FLM satisfies hypotheses {\bf H1} and {\bf H2'}. 
To check {\bf H1} one should 
check that the one dimensional Logistic Map intersects transversally 
$W^s(\Phi,\RR)$ the stable manifold of the renormalization operator. This is
an implicit assumption when one uses the renormalization
operator to explain the universality observed for the
Logistic Map. The only proof (to our knowledge) of this
fact is the one given by Lyubich (theorem 4.11 of \cite{Lyu99}). 
This proof is done in the space of quadratic-like germs,
which is a smaller space than the one considered here.
Hypothesis {\bf H2'} is trivial to check for 
the maps (\ref{FLM new set up}) and (\ref{FLM B new set up})).

On the other hand, note that theorem 3.8 in \cite{JRT11a} not only 
gives the existence of reducibility loss bifurcations, but it 
also gives an explicit expression of its slopes in term of the 
renormalization operator (the ones given by formulas 
(\ref{equation alpha n +}) and (\ref{equation alpha n -})). 
Actually we have given even a more explicit formula in 
corollary 3.13 in \cite{JRT11a}. We used these formulas 
to compute the reducibility directions $\alpha_n'(\omega)$ of 
the Forced Logistic Map (\ref{FLM new set up}). 

The initial values $\alpha_n$ have been computed numerically, 
by means of a Newton method applied to their invariance equation. 
Differentiating on formula (\ref{FLM new set up}) 
(respectively (\ref{FLM B new set up})) it 
is easy to write the values of $f^{(n)}_{0}, u^{(n)}_{0}$ and $v^{(n)}_{0}$ in 
terms of $\alpha_n$. Then using the discretization of the operator done in
Appendix \ref{section numerical computation of the specturm of L_omega},
we can compute numerically the iterates $f^{(n)}_{k}, u^{(n)}_{k}$ and $v^{(n)}_{k}$,
for $k=1,...,n-1$. Once we have these functions, we
can evaluate them to compute the values of $\alpha'_n(\omega) $
given by formula ($48$) of \cite{JRT11a}.

\begin{table}[t!]
\begin{center}
\texttt{\begin{tabular}{|c|c|c|c|}
\hline 
\rule{0pt}{2.5ex} n & $\alpha'_n(\omega)$  & $\bar{\epsilon}_a$ &
$\bar{\epsilon}_r $
\\ \hline
1 & -5.832915e+00 & 1.986668e-15 & 3.405962e-16 \\ 
2 & -8.494260e+00 & 1.155106e-13 & 1.359866e-14 \\ 
3 & -1.635128e+01 & 1.265556e-14 & 7.739800e-16 \\ 
4 & -1.125246e+01 & 3.995769e-14 & 3.551018e-15 \\ 
5 & -1.224333e+01 & 1.235427e-12 & 1.009061e-13 \\ 
6 & -1.807969e+01 & 8.989082e-12 & 4.971921e-13 \\ 
7 & -3.473523e+01 & 8.087996e-11 & 2.328470e-12 \\ 
8 & -2.958331e+01 & 2.204433e-10 & 7.451609e-12 \\ 
9 & -4.156946e+01 & 1.293339e-09 & 3.111273e-11 \\ 
10 & -7.896550e+01 & 3.847537e-08 & 4.872428e-10 \\ 
11 & -7.450073e+01 & 7.779131e-08 & 1.044168e-09 \\ 
\hline
\end{tabular}}
\end{center}
\vspace{-4mm}
\caption{
Values of $\alpha'_n(\omega)$ for the
FLM map (\ref{FLM new set up}) and $ \omega =  \frac{\sqrt{5}}{2}$. The
values  $\bar{\eps}_a$ and $\bar{\eps}_r$ correspond to the discrepancy with the
estimates of \cite{JRT11p2}  in absolute and relative terms.}
\label{taula DifDirMin}
\end{table}

We have used this procedure to compute the values of $\alpha'_n(\omega)$
for the map (\ref{FLM new set up}) (consequently also for the map
(\ref{q.p- family}) with $h_{\alpha,\eps}(\theta,x) = \alpha x (1-x) 
\cos(2\pi \theta)$) and  $\omega_0 = \frac{\sqrt{5}-1}{2}$.
The results are shown in table \ref{taula DifDirMin}. The values 
$ \alpha'_n(\omega)$ have been also computed
in \cite{JRT11p2} via a completely different 
procedure, based on a continuation method with extended precision. 
More concretely these values are displayed in table $2$ of the 
cited paper. The third and fourth columns of table \ref{taula DifDirMin}
display the discrepancies  between both computation, in
absolute ($\bar{\epsilon}_a$) and relative ($\bar{\epsilon}_r$)  terms.
This experiment has been repeated for other values of $\omega$ and also 
for the map (\ref{FLM B new set up}). In all cases we obtained that the 
slopes computed with both methods are the same up to similar accuracies 
to the ones displayed in table \ref{taula DifDirMin}.

This supports the correctness of both computations and at the same time
the conjecture {\bf \ref{conjecture H2}}, which has been assumed to be true 
to derive the formula used for this estimation.

\begin{appendix}

\section{Numerical computation of the spectrum of $\LL_\omega$}
\label{section numerical computation of the specturm of L_omega} 

In this appendix we present the numerical method that we used to 
discretize $\LL_\omega$ and to study its spectrum numerically. 
Our method is a slight modification of the one introduced by 
Lanford in \cite{Lan82} (see also \cite{Lan92}).

Let $\D(z_0, \rho )$ be the complex disc centered on 
$z_0$ with radius $\rho$. Consider $\RHH(\D(z_0, \rho ))$ the space 
of real analytic functions such that they are holomorphic on 
$\D(z_0, \rho )$ and continuous on it closure. Given a 
function $f\in \RHH(\D(z_0, \rho ))$, we can consider the following
modified Taylor expansion of $f$ around $z_0$, 
\begin{equation}
\label{modified taylor expansion} 
f(z) = \sum_{k=0}^{\infty} f_i \left(\frac{z-z_0}{\rho}\right)^i. 
\end{equation}

The truncation of a Taylor series at order $N$ induces 
a projection defined as
\[
\begin{array}{rccc}
p_{(N)}: &  \RHH(\D_{\rho}) & \rightarrow &  \R^{N+1} \\
   &       f          & \mapsto     &  (f_0, f_1, \dots,f_N) ,
\end{array}
\]

On the other hand we have its pseudo-inverse by the left
\[
\begin{array}{rccc}
i_{(N)}: &  \R^{N+1}           & \rightarrow & \RHH(\D_{\rho}) \\
       & (f_0, f_1, \dots,f_N)& \mapsto     & \displaystyle 
\sum_{k=0}^{N} f_i \left(\frac{z-z_0}{\rho}\right)^i ,
\end{array}
\]
in the sense that $i_{(N)} \circ p_{(N)} $ is the identity on $\R^{N+1}$. Note also
that both maps are linear.

Given a map $T:\RHH(\D(z_0, \rho ))\rightarrow \RHH(\D(z_0, \rho ))$, we can 
approximate it by its discretization 
$T^{(N)}: \R^{N+1} \rightarrow \R^{N+1}$ defined as
$T^{(N)}:= p_{(N)} \circ T \circ i_{(N)}$. More concretely if 
$T$ is a linear bounded operator,  we can 
compute the eigenvalues of the operator $T^{(N)}$ 
 in order to study the spectrum of $T$.  In general the
eigenvalues of $T^{(N)}$ might have nothing to do with the spectrum
of $T$. For example an infinite-dimensional operator does not need to
have eigenvalues, but a finite-dimensional one will always
have the same number of eigenvalues (counted with multiplicity)
as the dimension of the space. For this reason will do 
some numerical test on the results obtained with this discretization. 

At this point consider the map  $\LL_{\omega}: \RHH(\W)
\oplus\RHH(\W)  \rightarrow \RHH(\W) \oplus \RHH(\W)$ defined 
by equation (\ref{equation maps L_omega}). 
If we set $\W=\D(z_0,\rho)$ we can use 
the method described above in each component of 
$\RHH(\W) \oplus \RHH(\W)$ to discretize $\LL_\omega$ and 
approximate it by a map  
$\LL_{\omega}^{(N)} : \R^{2(N+1)} \rightarrow \R^{2(N+1)}$. 
Concretely in our computation we have taken $z_0=\frac{1}{5}$ and 
$\rho = \frac{3}{2}$. In figure \ref{figure inclusion boundaries} 
we include graphical evidence that the set $\W= \D\left( \frac{1}{5}, \frac{3}{2} 
\right)$ satisfies {\bf H0}, in section 
\ref{section definition and basic properties} can be found more details about 
this.

\begin{table}[t]
\begin{center}
\texttt{\begin{tabular}{|c|c||c|c|}
\hline 
$i$ &  $\lambda_i$  & $i$ & $\lambda_i$  	\\ \hline
1 & +7.8412640  +1.5617754$i$ & 13 & -0.0637772  +0.0000000$i$   \\
2 & +7.8412640  -1.5617754$i$  & 14 & -0.0637772  -0.0000000$i$   \\
3 & -2.5029079  +0.0000000$i$  & 15 & +0.0430641  +0.0435724$i$   \\
4 & -2.5029079  +0.0000000$i$  & 16 & +0.0430641  -0.0435724$i$   \\
5 & +0.5114250  +0.1942111$i$  & 17 & -0.0178305  +0.0165287$i$   \\
6 & +0.5114250  -0.1942111$i$  & 18 & -0.0178305  -0.0165287$i$   \\
7 & +0.4881230  +0.4930710$i$  & 19 & -0.0101807  +0.0000000$i$   \\
8 & +0.4881230  -0.4930710$i$  & 20 & -0.0101807  -0.0000000$i$   \\
9 & -0.3995353  +0.0000000$i$  & 21 & +0.0075181  +0.0069602$i$   \\
10 & -0.3995353  +0.0000000$i$  & 22 &  +0.0075181  -0.0069602$i$  \\ 
11 & -0.0982849  +0.0869398$i$  & 23 & -0.0029419  +0.0027336$i$  \\ 
12 & -0.0982849  -0.0869398$i$  & 24 & -0.0029419  -0.0027336$i$  \\ 
\hline
\end{tabular}}
\end{center}
\vspace{-4mm}
\caption{
The first twenty four eigenvalues of $\LL_{\omega}^{(N)}$, 
for $\omega=\frac{1-\sqrt{5}}{2}$. For the computation $N$ has 
been taken equal to 100.}
\label{table eigenvaules golden}
\end{table}

In table \ref{table eigenvaules golden} we have the first $24$ 
eigenvalues of $\LL_{\omega}^{(N)}$ for $N=100$ and
$\omega=\frac{\sqrt{5}-1}{2}$. The eigenvalues have been sorted by their 
modulus, from bigger to smaller. Note that the eigenvalues of
the discretized operator also satisfy  the properties 
given in section \ref{section The Fourier 
expansion of DT}. To justify the validity of these eigenvalues 
we have done the following numerical tests.

\begin{figure}[t!]
\begin{center}
\includegraphics[width=7.5cm]{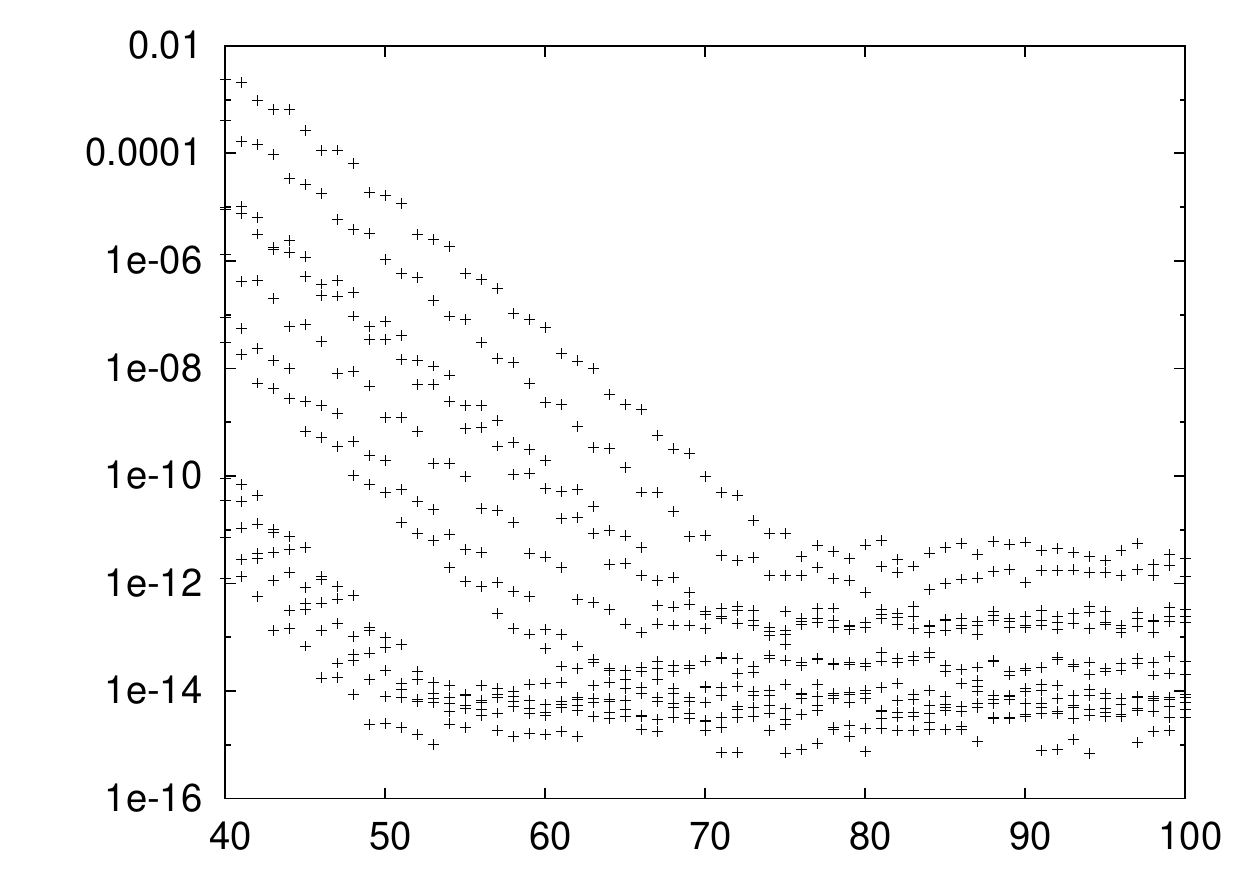}
\includegraphics[width=7.5cm]{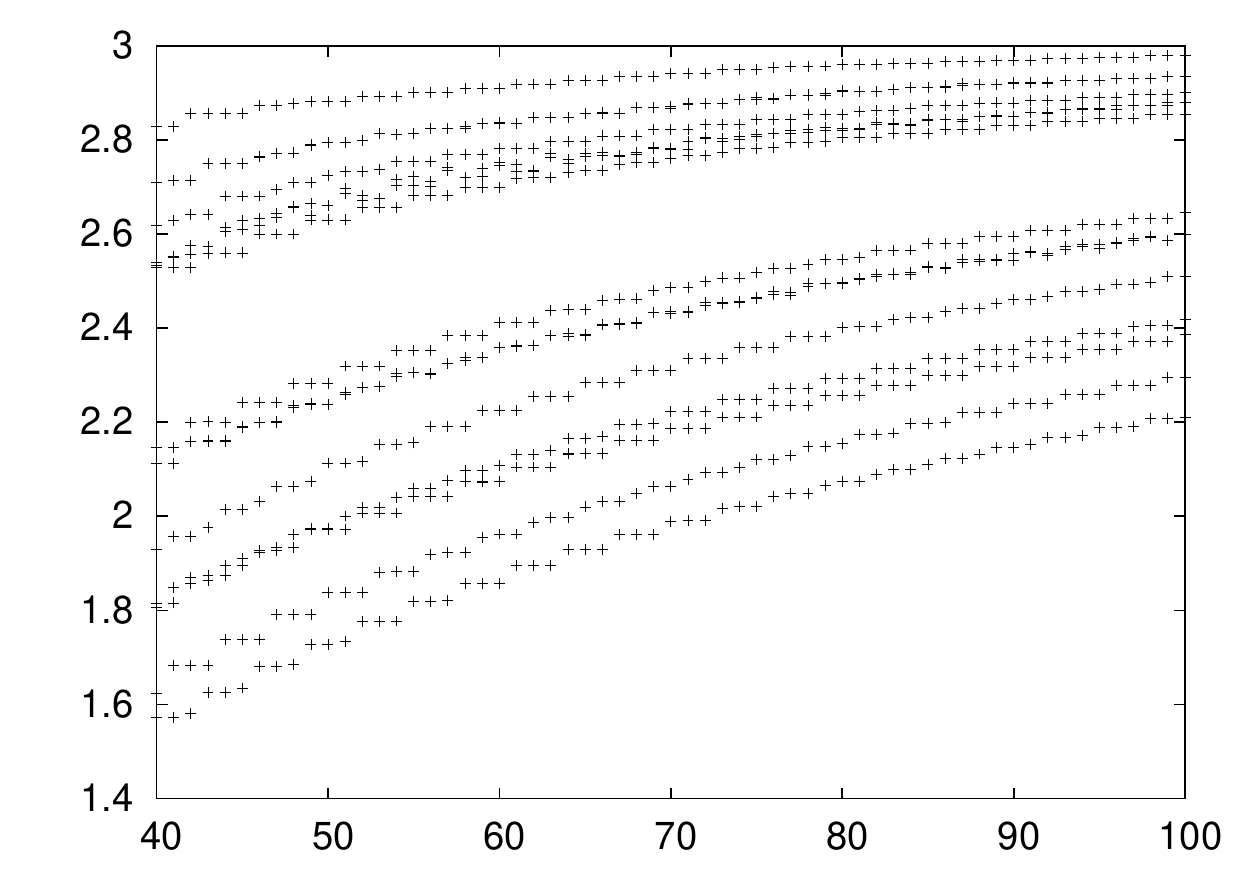}
\caption{
Estimation of the errors and the radii of convergence of 
the first twenty-four eigenvalues of $\LL_{\omega}$ for 
$\omega=\frac{\sqrt{5}-1}{2}$ with respect to the order 
of the discretization. See the text for more details. }
\label{Figure errors and radis golden}
\end{center}
\end{figure}

Consider that we have a real eigenvalue of multiplicity two, 
or a pair of complex eigenvalues\footnote{Note that there are no 
simple eigenvalues due to corollary 2.18 in \cite{JRT11a}}
which is persistent for different values of $N$ (the order of the 
discretization). The first test done to the eigenvalues is to 
check if the distance between the associated eigenvectors decreases 
when $N$ is increased.  In the left part of 
figure \ref{Figure errors and radis golden} we have 
the distance between the eigenvectors associated to the same 
eigenvalue of the operators $\LL^{(N)}_\omega$ and $\LL^{(110)}_\omega$ as 
a graph of $N$, with $N$ varying from $40$ to $100$. We have plotted 
this distance for the first twenty-four 
eigenvalues. To compute the distance between eigenvectors we have 
estimated the supremum norm of the difference between the real function 
represented by each of the vectors, in other words we have 
computed $\|i^{(N)}(v^{(N)}) - i^{(110)}(v^{(110)}) \|_\infty$ in 
the interval $(z_0-\rho,z_0+\rho)=\W\cap \R$. 

Note that, since the distance goes to zero this indicates that 
the eigenvectors, namely $v^{(N)}$,  converge to a 
limit $v^*$. One should expect these 
eigenvalues to be in the spectrum of $\LL_\omega$, 
but nothing ensures that $v^*$ belongs to the domain of $\LL_\omega$. 
We have done a second test on the reliability of the 
approximated eigenvectors, where we check this condition. 

Let us remark that with the numerical 
computations done so far, we have only checked that the eigenvectors
as elements of $  \RHH(\D(z_0,\rho))\oplus \RHH(\D(z_0,\rho)) $ to converge 
in the segment $(z_0-\rho,z_0+\rho)\subset \R$ instead of checking the 
convergence in the whole set $\D(z_0,\rho)$.  Let us give evidences that 
approximate eigenvectors obtained with our computations have a 
domain of analicity containing $\D(z_0,\rho)$. 

Consider that we have a function $f$ holomorphic in a domain of the complex 
plane containing $\D(z_0,\rho)$. If the we consider the expansion of 
$f$ given by equation (\ref{modified taylor expansion}), we have that 
$r$ the radius of convergence of the series around $z_0$ is given as
\[
r= \frac{\rho}{\limsup_{n \rightarrow \infty} \sqrt[n]{|f_n|} } . 
\] 
With the discretization considered here we have an approximation 
of the terms $f_n$, hence these can be used to compute a numerical 
estimation of $r$. 

Consider $v$ an eigenvector of the operator $\LL_\omega$. We have 
that $v=(v_1,v_2) \in \RHH(\W)\oplus \RHH(\W) $. Given 
$v_1^{(N)}= (v_1^{(N)}, v_2^{(N)})$ a numerical approximation of the 
eigenvector, we can use the procedure described above 
to estimate the radius of convergence of each  $v_1$ and $v_2$. 
We have done this for the eigenvectors associated to each of the first  twenty-four 
eigenvalues of  $\LL_\omega$ with $\omega= \frac{\sqrt{5}-1}{2}$ 
(keeping only the smaller of the two radius obtained). 
The results are displayed on the 
right part of figure \ref{Figure errors and radis golden}, 
where the estimated radius has been plotted 
with respect to $N$, the order of the discretization. Note that the estimations give a 
radius bigger than $\rho=\frac{3}{2}$, which indicates that the eigenvectors 
are analytic in $\D(z_0,\rho)$, and continuous on its closure, for 
$z_0=\frac{1}{5}$ and $\rho=\frac{3}{2}$.

\begin{figure}[t]
\begin{center}
\includegraphics[width=7.5cm]{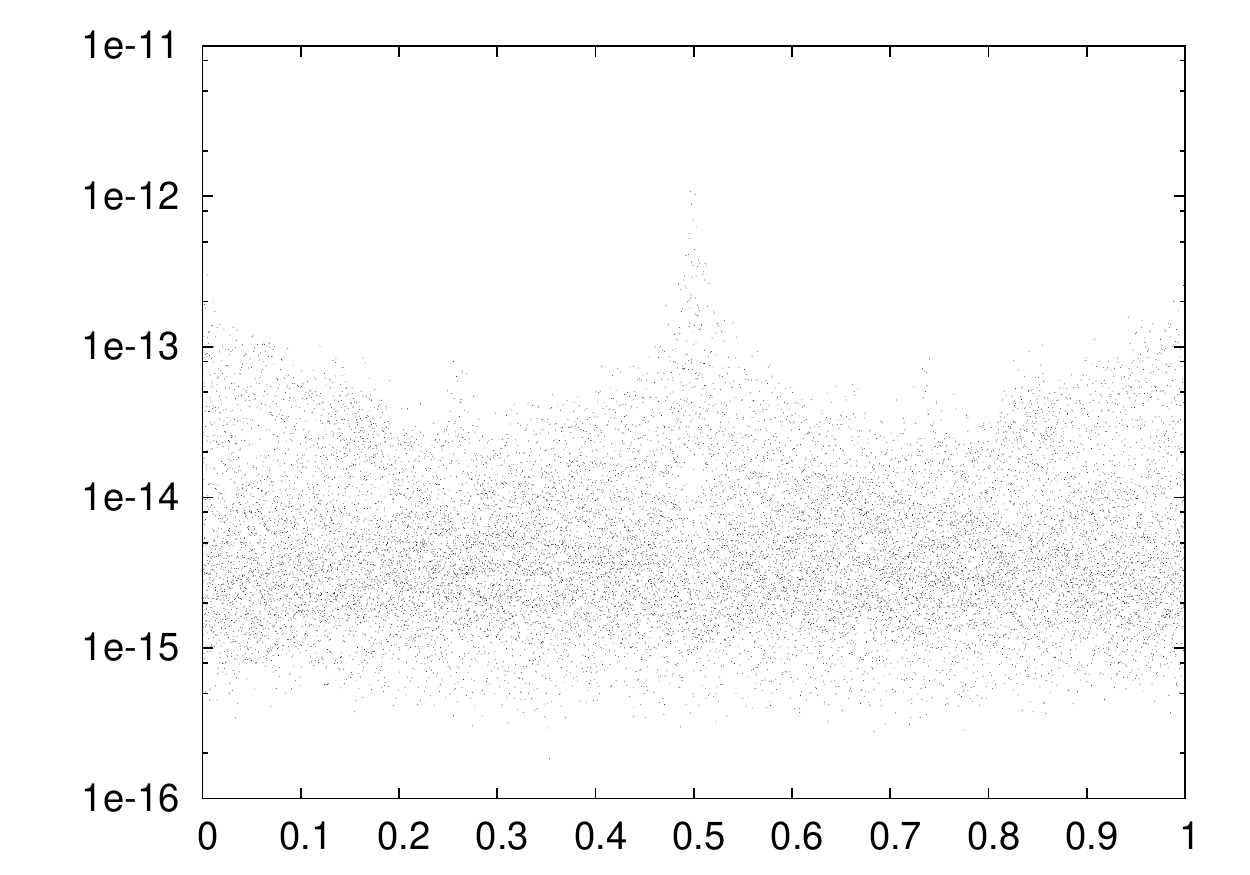}
\includegraphics[width=7.5cm]{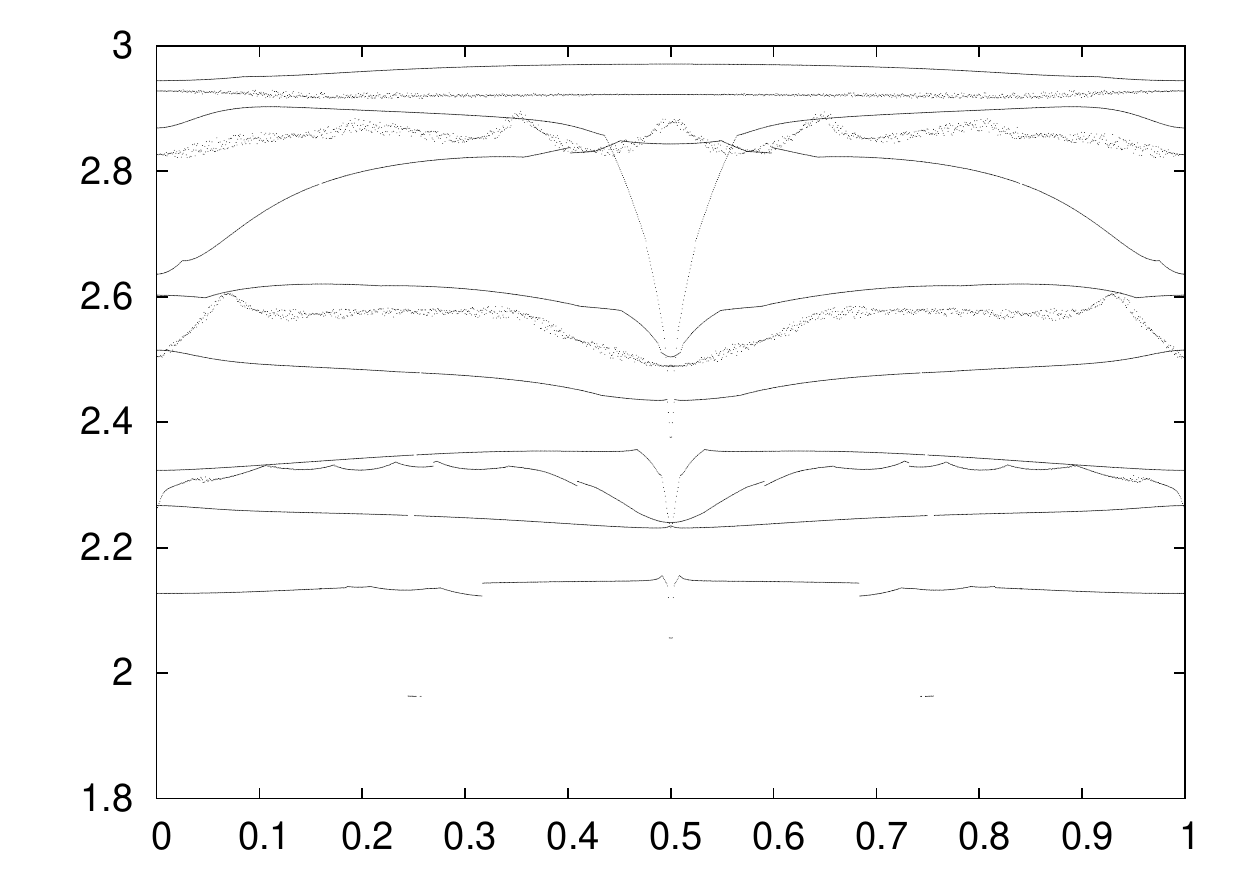}
\caption{
Estimation of the distance between eigenvectors with the 
same eigenvalue (left) 
and  estimation of the radios of convergence (right) of
the first twenty-four eigenvalues of $\LL_{\omega}$ with 
respect to $\omega$. See the text for more details. }
\label{Figure errors and radis omega}
\end{center}
\end{figure}

Up to this point, we have considered $\omega$ fixed to $\frac{\sqrt{5}-1}{2}$, 
but the same computations can be done to study the spectrum of 
$\LL_\omega$ with respect to the parameter $\omega$. In figure 
\ref{Operator Spectrum} of 
section \ref{section The Fourier expansion of DT} 
we have plotted  the first twenty-four eigenvalues 
of the map with respect to $\omega$. The set $\T$ has been discretized in a 
equispaced grid of $1280$ points.  Recall that the operator $\LL_{\omega}$ 
depends analytically on $\omega$ (proposition 2.20 in \cite{JRT11a}),  
therefore the spectrum also does (as long as the eigenvalues do not collide, 
see theorems III-6.17 and VII-1.7 in \cite{Kat66}). The numerical results 
agree with this analytic dependence. 

For this computation we have also made the same test as before to the eigenvalues.
The results of these tests are shown in figure \ref{Figure errors and radis omega}. 
To estimate the convergence of the eigenvectors we have 
compared the eigenspaces of the eigenvalues of $\LL_\omega^{(90)}$ with the 
eigenspaces associated to $\LL_\omega^{(100)}$ for each value of  $\omega$ in 
the cited grid of points on $\T$. The estimation of the radius of convergence has 
been also done with respect to $\omega$ for $N=90$. We 
have plotted the estimated error and convergence radius for the first
twenty-four eigenvalues in the same figure. Both result indicate 
that the eigenvalues obtained are reliable.

\end{appendix}

\bibliographystyle{plain}

\end{document}